\documentclass[a4paper,11pt]{article}
\usepackage{a4wide}
\usepackage[ansinew]{inputenc}
\usepackage{amssymb}
\usepackage{amsmath}
\usepackage{amsthm}
\usepackage{bbm}
\usepackage{mathrsfs}
\usepackage{enumerate}
\usepackage{needspace}

\usepackage[numbers, sort&compress]{natbib}
\input{math.sty}
\input{math2.sty}

\newcommand{\switch}{0}

\begin{document}

\title{A uniform central limit theorem and\\ efficiency for deconvolution estimators}

\author{Jakob S\"ohl\thanks{{The authors thank Richard Nickl and Markus Rei{\ss} for helpful comments and discussions. The first author thanks the Statistical Laboratory in Cambridge for the hospitality during a stay in which part of the project was started.
This research was supported by the Deutsche Forschungsgemeinschaft through the SFB 649 ``Economic Risk''.
}}~\thanks{E-mail address: soehl@math.hu-berlin.de} \and Mathias Trabs\footnotemark[1]~\thanks{E-mail address: trabs@math.hu-berlin.de}
}

\date{Humboldt-Universit\"at zu Berlin
\linebreak
\linebreak
\today}

\maketitle

\begin{abstract}
We estimate linear functionals in the classical deconvolution problem by kernel estimators. We obtain a uniform central limit theorem with $\sqrt{n}$--rate on the assumption that the smoothness of the functionals is larger than the ill--posedness of the problem, which is given by the polynomial decay rate of the characteristic function of the error. The limit distribution is a generalized Brownian bridge with a covariance structure that depends on the characteristic function of the error and on the functionals. The proposed estimators are optimal in the sense of semiparametric efficiency. The class of linear functionals is wide enough to incorporate the estimation of distribution functions. The proofs are based on smoothed empirical processes and mapping properties of the deconvolution operator. 
\end{abstract}


\section{Introduction}

Our observations are given by $n\in\N$ independent and identically distributed random variables
\begin{equation}\label{eqObservation}
  Y_j=X_j+\eps_j, \quad \quad j=1,\dots,n,
\end{equation}
where $X_j$ and $\eps_j$ are independent of each other, the distribution of the errors $\eps_j$ is supposed to be known and the aim is statistical inference on the distribution of $X_j$. Let us denote the densities of $X_j$ and $\eps_j$ by $f_X$ and $f_\eps$, respectively. We consider the case of ordinary smooth errors, which means that the characteristic function $\phi_\eps$ of the errors $\eps_j$ decays with polynomial rate, determining the ill--posedness of the inverse problem.
The contribution of this article to the well studied problem of deconvolution is twofold. First, we prove a uniform central limit theorem for kernel estimators of the distribution function of $X_j$ in the setting of $\sqrt n$ convergence rates.
More precisely, the theorem does not only include the estimation of the distribution function but covers translation classes of linear functionals of the density $f_X$ whenever the ill--posedness is smaller than the smoothness of the functionals.
Second, we obtain more exact results than the minimax rates of convergence by showing that the used estimators are optimal in the sense of semiparametric efficiency.

The classical Donsker theorem plays a central role in statistics and states that the empirical distribution function of an independent, identically distributed sample converges uniformly to the distribution function. In the deconvolution model \eqref{eqObservation} our Donsker theorem states  uniform convergence for an asymptotically unbiased estimator of translated functionals $t \mapsto \theta_t:=\int\zeta(x-t) f_X(x)\d x$, where the special case $\zeta:=\mathbbm{1}_{(-\infty,0]}$ leads to the estimation of the distribution function. This generalization allows to consider functionals $\theta_t$ as long as the smoothness of $\zeta$ in an $L^2$--Sobolev sense compensates the ill--posedness of the problem.
The limiting process $\mathbb G$ in the uniform central limit theorem is a generalized Brownian bridge, whose covariance depends on the functional $\zeta$ and through the deconvolution operator $\F^{-1}[1/\phi_\eps]$ also on the distribution of the errors.
The used kernel estimators $\widehat \theta_t$ are minimax optimal since they converge with a $\sqrt{n}$--rate. So investigating optimality further leads naturally to the question whether the asymptotic variance of the estimators is minimal, as in the case of the empirical distribution function in the classical Donsker theorem. We prove that the estimator $\widehat \theta_{\bull}$ is efficient in the sense of a H\'ajek--Le Cam convolution theorem. In particular, the asymptotic covariance matrices of the finite dimensional distributions achieve the Cram\'er--Rao information bound. By uniform convergence and efficiency the kernel estimator of $f_X$ fulfills the `plug-in' property of \citet{BickelRitov2003} in the deconvolution model \eqref{eqObservation}.

The deconvolution problem has attracted much attention so we mention here only closely related works and refer the interested reader to the references therein. The classical works by \citet{fan1991b,fan1991} contain asymptotic normality of kernel density estimators as well as minimax convergence rates for estimating the density and the distribution function. \citet{ButuceaComte2009} have treated the data--driven choice of the bandwidth for estimating functionals of $f_X$ but assumed some minimal smoothness and integrability conditions on the functional $\theta_t$, which exclude, for example, $\zeta:=\mathbbm{1}_{(-\infty,0]}$ since it is not integrable.
\citet{DattnerEtAll2011} have studied minimax--optimal and adaptive estimation of the distribution function. Asymptotic normality of estimators for the distribution function has been shown by \citet{EsUh2005} in the case of supersmooth errors an by \citet{HallLahiri2008} for ordinary smooth errors. In contrast we consider the estimation of general linear functionals and are interested in uniform convergence. Uniform results have been studied for the density but not for the distribution function by \citet{BissantzEtAl2007} and  by \citet{LouniciNickl2011}. Recently, \citet{NicklReiss2012} have proved a Donsker theorem for estimators of the distribution function of a L\'evy measure. Their situation is related but more involved than ours, owing to the nonlinearity and the auto-deconvolution of the L\'evy measure. In a deconvolution context we consider the more general problem of estimating linear functionals efficiently, which contains estimating of the distribution function as a special case and provides clear insight in the interplay between smoothness of~$\zeta$ and the ill--posedness of the problem. While efficiency has been  investigated in various semiparametric models, e.g., see \citet{bicklEtAl1998}, to the best of the authors knowledge there are no results in this direction in the deconvolution framework. However, in the L\'evy setting \citet{NicklReiss2012} have shown heuristically that their estimator achieves the lower bound of the variance while a rigorous proof remained open.

In order to show the uniform central limit theorem in the deconvolution problem, we prove that the empirical process $\sqrt{n}(\PP_n-\PP)$ is tight in the space of bounded functions acting on the class
\[
\G:=\{\F^{-1}[1/\phi_\eps(-\bull)]\ast\zeta_t| \: t\in\R\},\qquad\zeta_t:=\zeta(\bull-t),
\]
where $\PP$ and $\PP_n=\tfrac{1}{n}\sum_{j=1}^{n}\delta_{Y_j}$ denote the true and the empirical probability measure of the observations $Y_j$, respectively.
Since $\G$ may consist of translates of an unbounded function, this is in general not a Donsker class. Nevertheless, \citet{RadulovicWegkamp2000} have observed that a smoothed empirical processes might converge even when the unsmoothed process does not. \citet{GineNickl2008} have further developed these ideas and have shown uniform central limit theorems for kernel density estimators. \citet{NicklReiss2012} used smoothed empirical processes in the inverse problem of estimating the distribution function of L\'evy measures.
In order to show semiparametric efficiency in the deconvolution problem, the main problem is to show that the efficient influence function is indeed an element of the tangent space. If the regularity of $\zeta$ is small, the standard methods given in the monograph of \citet{bicklEtAl1998} do not apply in this ill--posed problem. Instead, we approximate $\zeta$ by a sequence of smooth $(\zeta_n)$ and show the convergence of the information bounds. Interestingly, this reveals a relation between the intrinsic metric of the limit $\mathbb G$ and the metric which is induced by the inverse Fisher information.
Additionally to techniques of smoothed empirical processes and the calculus of information bounds, our proofs rely on the Fourier multiplier property of the underlying deconvolution operator $\F^{-1}[1/\phi_\eps]$, which is related to pseudo-differential operators as noted in the L\'evy process setting by \citet{NicklReiss2012} and in the deconvolution context by \citet{SchmidtHieberEtAll2012}. Important for our proofs are the mapping properties of $\F^{-1}[1/\phi_\eps]$ on Besov spaces.\par
This paper is organized as follows: In Section~\ref{secModel} we formulate the Donsker theorem and discuss its consequences. Efficiency is then considered in Section~\ref{secEfficiency}. All proofs are deferred to Sections~\ref{secProofs} and \ref{secProofs2}. In the Appendix we summarize definitions and properties of the function spaces used in the paper. 
\section{Uniform central limit theorem}\label{secModel}

\subsection{The estimator}
According to the observation scheme \eqref{eqObservation}, $Y_j$ are distributed with density $ f_Y=f_X\ast f_\eps$ determining the probability measure $\PP$. The characteristic function~$\phi$ of~$\PP$ can be estimated by its empirical version $\phi_n(u)=\frac{1}{n}\sum_{j=1}^ne^{iuY_j}, u\in\R$. For $\zeta$ to be specified later and recalling $\zeta_t=\zeta(\bull-t)$, our aim is to estimate functionals of the form
\begin{equation}
  \theta_t:=\langle \zeta_t,f_X\rangle=\int \zeta_t(x) f_X(x)\d x.
\end{equation}
Defining the Fourier transform by $\F f(u):=\int e^{iux}f(x)\d x, u\in\R$, the natural estimator of the functional $\theta_t$ is given by
\begin{equation}\label{eqEstimator}
  \widehat \theta_t:=\int\zeta_t(x)\F^{-1}\Big[\F K_h\frac{\phi_n}{\phi_\eps}\Big](x)\d x,
\end{equation}
where $K$ is a kernel, $h>0$ the bandwidth and we have written as usual $K_h(x)=h^{-1}K(x/h)$. Choosing $\F K=\mathbbm1_{[-\pi,\pi]}$ for some $\pi>0$ leads to the estimator proposed by \citet{ButuceaComte2009}. Throughout, we suppose that
\begin{enumerate}
  \item $K\in L^1(\R)\cap L^\infty(\R)$ is symmetric and band--limited with $\supp(\F K)\subset[-1,1]$,
  \item for $l=1,\dots, L$
  \begin{equation}\label{eqKorder}
    \int K=1,\quad\int x^lK(x)\d x=0,\quad\int |x^{L+1}K(x)|\d x<\infty\quad\text{and}
  \end{equation}
  \item $K\in C^1(\R)$ satisfies, denoting $\langle x\rangle:=(1+x^2)^{1/2}$,
    \begin{equation}\label{assKDecay}
      |K(x)|+|K'(x)|\lesssim\langle x\rangle^{-2}.
    \end{equation}
\end{enumerate}
Throughout, we write $A_p\lesssim B_p$ if there exists a constant $C>0$ independent of the parameter $p$ such that $A_p\le CB_p$. If $A_p\lesssim B_p$ and $B_p\lesssim A_p$, we write $A_p\sim B_p$. Examples of such kernels can be obtained by taking $\F K$ to be a symmetric function in $C^\infty(\R)$ which is supported in $[-1,1]$ and constant to one in a neighborhood of zero. The resulting kernels are called flat top kernels and were used in deconvolution problems, for example, by \citet{BissantzEtAl2007}.

\subsection{Statement of the theorem}
Given a function $\zeta$ specified later, our aim is to show a Donsker theorem for the estimator over the class of translations $\zeta_t$, $t\in \R$. In view of the classical Donsker theorem in a model without additive errors, where no assumptions on the smoothness of the distribution are needed, we want to assume as less smoothness of $f_X$ as possible still guaranteeing $\sqrt n$-rates. For some $\delta>0$ the following assumptions on the density $f_X$ will be needed:
\begin{assumption}\label{assFx}\needspace{5\baselineskip}
  \hspace{1em}
  \begin{enumerate}
    \item\label{assFxMoment} Let $f_X$ be bounded and assume the moment condition $\int |x|^{2+\delta}f_X(x)\d x<\infty$.
    \item\label{assFxSmooth} Assume $f_X\in H^\alpha(\R)$ that is the density has Sobolev smoothness of order $\alpha\ge0$.
  \end{enumerate}
\end{assumption}
We refer to the appendix for an exact definition of the Sobolev space $H^\alpha(\R)$. Boundedness of the observation density $f_Y$ follows immediately from~\ref{assFxMoment} since $\|f_Y\|_\infty\le\|f_X\|_\infty\|f_\eps\|_{L^1}<\infty$. In addition to the smoothness of $f_X$, the smoothness of $\zeta$ will be crucial. We assume for $\gamma_s,\gamma_c>0$
\begin{align}
  \zeta\in Z^{\gamma_s,\gamma_c}:=\Big\{\zeta=&\zeta^c+\zeta^s \Big| \zeta^s\in H^{\gamma_s}(\R) \text{ is compactly supported as well }\notag\\
    &\text{as }\langle x \rangle ^\tau\big(\zeta^c(x)-a(x)\big)\in H^{\gamma_c}(\R)\label{eqZ}\text{ for some }\tau>0\text{ and }\\
    &\text{some }a\in C^\infty(\R)\text{ such that }a'\text{ is compactly supported}\Big\}\notag
\end{align}
and write for $\zeta \in Z^{\gamma_s,\gamma_c}$ with a given decomposition $\zeta=\zeta^s+\zeta^c$
\[
  \|\zeta\|_{Z^{\gamma_s,\gamma_c}}:=\|\zeta^s\|_{H^{\gamma_s}}+\big\|\tfrac{1}{ix+1}\zeta^c(x)\big\|_{H^{\gamma_c}},
\]
\ifodd\switch
  which is finite since the term $\|\tfrac{1}{ix+1}\zeta^c(x)\|_{H^{\gamma_c}}$ can be bounded by $\|\tfrac{a(x)}{ix+1}\|_{H^{\gamma_c}}+\|\tfrac{1}{(ix+1)\langle x\rangle^{\tau}}\|_{C^s}\|\langle x\rangle^{\tau}(\zeta^c(x)-a(x))\|_{H^{\gamma_c}}<\infty$ for any $s>\gamma_c$.
\else
  which is finite since $\|\tfrac{1}{ix+1}\zeta^c(x)\|_{H^{\gamma_c}}$ is bounded by $\|\tfrac{a(x)}{ix+1}\|_{H^{\gamma_c}}+\|\tfrac{1}{(ix+1)\langle x\rangle^{\tau}}\|_{C^s}\|\langle x\rangle^{\tau}(\zeta^c(x)-a(x))\|_{H^{\gamma_c}}<\infty$ for any $s>\gamma_c$.
\fi
Several examples for $\zeta$ and corresponding $\gamma_s,\gamma_c$ will be given in Examples 1-3 below. In particular, $\mathbbm 1_{(-\infty,0]}\in Z^{\gamma_s,\gamma_c}$ for $\gamma_s<1/2$. The ill--posedness of the problem is determined by the decay of the characteristic function of the errors. More precisely, we suppose
\begin{assumption}\label{assEps}\needspace{5\baselineskip}
  Let the error distribution satisfy
  \begin{enumerate}
    \item\label{assEpsMoment} $\int |x|^{2+\delta} f_\eps(x)\d x<\infty$ thus $\phi_\eps$ is twice continuously differentiable and
    \item\label{assEpsDecay} $|(\phi_\eps^{-1})'(u)|\lesssim\langle u\rangle^{\beta-1}$ for some $\beta>0$, in particular $|\phi_\eps^{-1}(u)|\lesssim\langle u\rangle^{\beta}, u\in\R$.
  \end{enumerate}
\end{assumption}
\noindent
Throughout, we write $\phi_\eps^{-1}=1/\phi_\eps$. The Assumption \ref{assEpsDecay} on the distribution of the errors is similar to the classical decay assumption by \citet{fan1991b} and it is fulfilled for many ordinary smooth error laws such as gamma or Laplace distributions as discussed below. Assumption~\ref{assEps}\ref{assEpsDecay} implies that $\phi_\eps^{-1}$ is a Fourier multiplier on Besov spaces so that
\[
  B_{p,q}^s(\R)\ni f\mapsto\F^{-1}[\phi_\eps^{-1}(-\bull)\F f]\in B_{p,q}^{s-\beta}(\R)
\]
for $p,q\in[1,\infty], s\in\R$, is a continuous linear map, which is essential in our proofs, compare Lemma~\ref{lemPsiDO}. In the same spirit \citet{SchmidtHieberEtAll2012} discuss the behavior of the deconvolution operator as pseudo--differential operator. We define
\begin{equation}\label{g}
g_t:=\F^{-1}[\phi^{-1}_\eps(-\bull)]\ast\zeta_t\quad\text{and}\quad\G=\{g_t|t\in\R\}.
\end{equation}
Note that in general $g_t$ may only exist in a distributional sense, but on Assumption~\ref{assEps} and for $\zeta\in Z^{\gamma_s,\gamma_c}$ it can be rigorously interpreted by (see \eqref{eqDecomp})
\begin{align*}
  g_0(x)=&\F^{-1}[\varphi_\eps^{-1}(-u)\F \zeta^s(u)](x)\notag\\
  &+(1+ix)\F^{-1}[\varphi_\eps^{-1}(-u)\F[\tfrac{1}{iy+1}\zeta^c(y)](u)](x)\notag\\
  &+\F^{-1}[(\varphi_\eps^{-1})'(-u)\F[\tfrac{1}{iy+1}\zeta^c(y)](u)](x),
\end{align*}
which indicates why we have imposed an assumption on $(\phi_\eps^{-1})'$ and have defined $\|\bull\|_{Z^{\gamma_s,\gamma_c}}$ as above.\par
It will turn out that $\G$ is $\PP$--pregaussian, but not Donsker in general.
Denoting by $\lfloor\alpha\rfloor$ the largest integer smaller or equal to $\alpha$ and
defining convergence in law on $\ell^\infty(\R)$ as \citet[p. 94]{dudley1999}, we state our main result
\begin{theorem}\label{thmDonsker}
  Grant Assumptions~\ref{assFx} and \ref{assEps} as well as $\zeta\in Z^{\gamma_s,\gamma_c}$ with $\gamma_s>\beta, \gamma_c>(1/2\vee\alpha)+\gamma_s$ and $\alpha+3\gamma_s>2\beta+1$. Furthermore, let the kernel $K$ satisfy \eqref{eqKorder} with $L=\lfloor\alpha+\gamma_s\rfloor$.
  Let $h_n^{2\alpha+2\gamma_s} n\to 0$ and if $\gamma_s\le\beta+1/2$ let in addition $h_n^\rho n\to\infty$ for some $\rho>4\beta-4\gamma_s+2$, then
  \[
    \sqrt n(\widehat\theta_t-\theta_t)_{t\in\R}\overset{\mathcal L}{\longrightarrow}\mathbb G \quad\text{in }\ell^\infty(\R)
  \]
  as $n\to\infty$, where $\mathbb G$ is a centered Gaussian Borel random variable in $\ell^\infty(\R)$ with covariance function given by
  \[
    \Sigma_{s,t}:=\int g_s(x)g_t(x)\PP(\d x)-\theta_s\theta_t
  \]
   for $g_s,g_t$ defined in \eqref{g} and $s,t\in\R$.
\end{theorem}
We illustrate the range of this theorem by the following examples.
\begin{examples}\label{exInd}
  We consider the indicator function $\mathbbm 1_{(-\infty,0]}(x)$, $x\in\R$. Let $a$ be a monotone decreasing $C^\infty(\R)$ function, which is for some $M>0$ equal to zero for all $x\ge M$ and equal to one for all $x\le -M$. We define $\zeta^s:=\mathbbm 1_{(-\infty,0]}-a$ and $\zeta^c:=a$. From the bounded variation of $\zeta^s$ follows $\zeta^s\in B^1_{1,\infty}(\R)\subset H^{\gamma_s}(\R)$ for any $\gamma_s<1/2$ by Besov smoothness of bounded variation functions \eqref{eqBoundVar} as well as by the Besov space embeddings \eqref{eqBpEmbed} and \eqref{eqBqEmbed}. Since $a\in C^\infty(\R)$ and $a'$ is compactly supported, the condition on $\zeta^c$ is satisfied for any $\gamma_c>0$. Hence, $\mathbbm 1_{(-\infty,t]}\in Z^{\gamma_s,\gamma_c}$ if $\gamma_s<1/2$. On the other hand, this cannot hold for $\gamma_s>1/2$ since $H^{\gamma_s}(\R)\subset C^0(\R)$ by Sobolev's embedding theorem or by \eqref{eqCsEmbed}, \eqref{eqBpEmbed} and \eqref{eqBqEmbed}. Owing to the condition $\gamma_s>\beta$, Assumption~\ref{assEps} needs to be fulfilled for some $\beta<1/2$ which is done, for example, by the gamma distribution $\Gamma(\beta,\eta)$ with $\beta\in(0,1/2)$ and $\eta\in(0,\infty)$, that is
  \[
      f_\eps(x):=\gamma_{\beta,\eta}(x):=\frac{1}{\Gamma(\beta)\eta^\beta}x^{\beta-1}e^{-x/\eta}\mathbbm 1_{[0,\infty)}(x),\quad x\in \R,
   \]
   and $\varphi_\eps(u)=(1-i\eta u)^{-\beta}$, $u\in\R$.
\end{examples}
\begin{examples}
   Let $\zeta_t(x):=\zeta^s_t(x):=\max(K-|x-t|,0)$ and $\zeta^c_t(x):=0$ with $K>0$. The payoff of the butterfly spread is described by such a function \citep{FoellmerSchied2004}. Then $\F\zeta(u)=4\sin^2(u/2)/u^2$ and $\zeta^s\in H^{\gamma_s}(\R)$ for any $\gamma_s<3/2$.  So, Assumption~\ref{assEps} is required for some $\beta<3/2$, which holds, for example, for the chi--squared distribution with one or two degrees of freedom or for the exponential distribution.
\end{examples}
\begin{examples}
  \citet{ButuceaComte2009} studied the case $\beta>1$ and derived $\sqrt n$-rates for $\gamma_s>\beta$ in our notation. In particular, they considered supersmooth~$\zeta$, that is $\F\zeta$ decays exponentially. In this case $\zeta\in H^s(\R)$ for any $s\in\N$. Requiring the slightly stronger assumption that $\langle x \rangle^\tau \zeta(x) \in H^s(\R)$ for some arbitrary small $\tau>0$ and for all $s\in \N$ we can choose $\zeta^c:=\zeta$ and $\zeta^s:=0$. Then~$\beta$ can be taken arbitrary large such that all gamma distributions, the Laplace distributions and convolutions of them can be chosen as error distributions.
\end{examples}

\subsection{Discussion}

To have $\sqrt n$--rates we suppose $\gamma_s>\beta$ in Theorem~\ref{thmDonsker}, which means that the smoothness of the functionals compensates the ill--posedness of the problem. This condition is natural in view of the abstract analysis in terms of Hilbert scales by \citet{GoldenshlugerPereverzev2003}, who obtain the minimax rate $n^{-(\alpha+\gamma_s)/(2\alpha+2\beta)}\vee n^{-1/2}$ in our notation. As a consequence of the condition on $\gamma_s$ and $\gamma_c$ we can bound the stochastic error term of the estimator $\widehat\theta_t$ uniformly in $h\in(0,1)$. The bias term is of order $h^{\alpha+\gamma_s}$.

For $\gamma_s>\beta+1/2$  the class $\mathcal G$ is a Donsker class. In this case the only condition on the bandwidth is that the bias tends faster than $n^{-1/2}$ to zero. In the interesting but involved case $\gamma_s\in(\beta,\beta+1/2]$, the class $\mathcal G$ will in general not be a Donsker class. Estimating the distribution function as in Example~\ref{exInd} belongs to this case. In order to see that $\mathcal G$ is in general not a Donsker class, let the error distribution be given by $f_\eps=\gamma_{\beta,\eta}(-\bull)$ and $\zeta=\gamma_{\sigma,\eta}$ with $\sigma\in(\gamma_s+1/2,\beta+1)$. Then $g_t$ equals $\gamma_{\sigma-\beta,\eta}\ast\delta_t$. For the shape parameter holds $\sigma-\beta\in(1/2,1)$ and thus $g_t$ is an $L^2(\R)$--function unbounded at $t$. The Lebesgue density of $\PP$ is bounded by Assumption~\ref{assFx}\ref{assFxMoment}. Hence, $\mathcal G$ consists of all translates of an unbounded function and thus cannot be Donsker, cf. Theorem~7 by \citet{nickl2006}.

Therefore, for $\gamma_s\in(\beta,\beta+1/2]$ smoothed empirical processes are necessary, especially we need to ensure enough smoothing to be able to obtain a uniform central limit theorem. The bandwidth cannot tend too fast to zero, more precisely we require $h_n^{\rho}n\to\infty$ as $n\to\infty$ for some $\rho$ with $\rho>4\beta-4\gamma_s+2$. In combination with the bias condition $h_n^{2\alpha+2\gamma_s}n\to0$ as $n\to\infty$ we obtain necessarily  $\alpha+\gamma_s>2\beta-2\gamma_s+1$ leading to the assumption in the theorem.
Since $2\alpha+2\gamma_s>\alpha+2\beta-\gamma_s+1>4\beta-4\gamma_s+2$ we can always choose $h_n\sim n^{-1/(\alpha+2\beta-\gamma_s+1)}$.
In contrast to \citet{fan1991, ButuceaComte2009, DattnerEtAll2011} our choice of the bandwidth $h_n$ is not determined by the bias--variance trade--off, but rather by the amount of smoothing necessary to obtain a uniform central limit theorem.
The classical bandwidth $h_n\sim n^{-1/(2\alpha+2\beta)}$ is optimal for estimating the density in the sense that it achieves the minimax rate with respect to the mean integrated squared error (MISE), compare \citet{fan1991} who assumes H\"older smoothness of $f_X$ instead of $L^2$--Sobolev smoothness. For this choice the bias condition $h_n^{2\alpha+2\gamma_s}n\to0$ is satisfied. If $\gamma_s\le\beta+1/2$ the classical bandwidth satisfies the additional minimal smoothness condition in the case of estimating the distribution function with mild conditions on $f_X$. It suffices for example that $f_X$ is of bounded variation. Then $\alpha$ and $\gamma_s$ can be chosen large enough in $(0,1/2)$ such that $2\alpha+2\beta>4\beta-4\gamma_s+2$ and the classical bandwidth satisfies the conditions of the theorem. Whenever the classical bandwidth $h_n\sim n^{-1/(2\alpha+2\beta)}$ satisfies the conditions of Theorem~\ref{thmDonsker}, then the corresponding density estimator is a `plug--in' estimator in the sense of \citet{BickelRitov2003} meaning that the density is estimated rate optimal for the MISE, the functionals are estimated efficiently (see Section~\ref{secEfficiency}) and the estimators of the functionals converge uniformly over $t\in\R$.

The smoothness condition on the density $f_X$ is then a consequence of the given choice of $h_n$ together with the classical bias estimate for kernel estimators. As we have seen in Example~\ref{exInd} for estimating the distribution function we have $\zeta=\mathbbm 1_{(-\infty,0]}\in Z^{\gamma_s,\gamma_c}$ with $\gamma_s<1/2$ arbitrary close to 1/2. In the classical Donsker theorem which corresponds to the case $\beta\to0$ the condition $\alpha+3\gamma_s>2\beta+1$ would simplify to $\alpha>-1/2$. However, we suppose $f_X$ to be bounded, which leads to much clearer proofs, and thus $f_X\in H^0(\R)$ is automatically satisfied. Assumption~\ref{assFx} allows to focus on the interplay between the functional $\zeta$ and the deconvolution operator $\F^{-1}[\phi_\eps^{-1}]$. \citet{NicklReiss2012} have studied the case of unbounded densities, which is necessary in the L\'evy process setup, but considered $\zeta_t=\mathbbm 1_{(-\infty,t]}$ only. The class $Z^{\gamma_s,\gamma_c}$ is defined by $L^2$--Sobolev conditions so that bounded variation arguments for $\zeta$ have to be avoided in the proofs.

An interesting aspect is the following: If we restrict the uniform convergence to $(\zeta_t)_{t\in T}$ for some compact set $T\subset\R$, it is sufficient to assume $\frac{1}{ix+1}\zeta^c\in H^{\gamma_c}(\R)$ instead of requiring $(1\vee |x|^\tau)(\zeta^c(x)-a(x))\in H^{\gamma_c}(\R)$ for some $\tau>0$ and a function $a\in C^\infty(\R)$ such that $a'$ is compactly supported as done in $Z^{\gamma_s,\gamma_c}$. In particular, slowly growing $\zeta$ would be allowed.
The stronger condition in the definition of $Z^{\gamma_s,\gamma_c}$ is only needed to ensure polynomial covering numbers of $\{g_t|t\in T\}$ for $T\subset\R$ unbounded (cf. Theorem~\ref{thmPregauss} below).\par

As a corollary of Theorem~\ref{thmDonsker} we can weaken Assumption~\ref{assEps}\ref{assEpsDecay}. If the characteristic function of the errors $\eps$ is given by $\tilde\phi_\eps=\phi_\eps\psi$ where $\phi_\eps$ satisfies Assumption~\ref{assEps}\ref{assEpsDecay} and there is a Schwartz distribution $\nu\in\mathscr S'(\R)$ such that $\F\nu=\psi^{-1}$ and $\nu\ast\zeta\in Z^{\gamma_s,\gamma_c}$ for $\zeta\in Z^{\gamma_s,\gamma_c}$, then for $t\in\R$
\[
  \F^{-1}[\tilde\phi_\eps^{-1}]\ast\zeta(\bull-t)=\F^{-1}[\phi_\eps^{-1}]\ast(\nu\ast\zeta)(\bull-t)
\]
and thus we can proceed as before. For instance, for translated errors $f_\eps\ast\delta_\mu$ with $\mu\neq0$, the distribution $\nu$ would be given by $\delta_{-\mu}$.

As for the classical Donsker theorem the Donsker theorem for deconvolution estimators has many different applications, the most obvious being the construction of confidence bands. Further Donsker theorems may be obtained by applying the functional delta method to Hadamard differentiable maps. Let us illustrate the construction of confidence bands.
By the continuous mapping theorem we infer
\begin{align*}
  \sup_{t\in\R}\sqrt{n}|\widehat\theta_t-\theta_t|\overset{\mathcal L}{\longrightarrow}\sup_{t\in\R} |\mathbb G(t)|.
\end{align*}
The construction of confidence bands reduces now to knowledge about the distribution of the
supremum of $\mathbb G$. Suprema of Gaussian processes are well studied and information about their distribution can be either obtained from theoretical considerations as in \citet[App. A.2]{vanderVaartWellner1996} or from Monte Carlo simulations. Let $q_{1-\alpha}$ be the $(1-\alpha)$--quantile of $\sup_{t\in\R} |\mathbb G(t)|$ that is
$\PP(\sup_{t\in\R} |\mathbb G(t)|\le q_{1-\alpha})= 1-\alpha$. Then
\begin{align*}
  \lim_{n\to\infty}\PP\left(\theta_t\in[\widehat\theta_t-q_{1-\alpha} n^{-1/2},\widehat\theta_t+q_{1-\alpha} n^{-1/2}] \text{ for all }t\in\R\right)= 1-\alpha
\end{align*}
and thus the intervals $[\widehat\theta_t-q_{1-\alpha} n^{-1/2},\widehat\theta_t+q_{1-\alpha} n^{-1/2}]$ define a confidence band.

\section{Efficiency}\label{secEfficiency}
Having established the asymptotic normality of our estimator, the natural question is whether it is optimal in the sense of the convolution Theorem 5.2.1 by \citet{bicklEtAl1998}. Typically, efficiency is investigated for estimators $T_n$ which are (locally) regular, that is for any  parametric submodel $\eta\to f_{X,\eta}$ and $n^{1/2}|\eta_n-\eta|\lesssim1$ the law of $n^{1/2}(T_n-\langle\zeta, f_{X,\eta}\rangle)$ under $\eta_n$ converges for $n\to\infty$ to a distribution independent of $(\eta_n)$. In Lemma~\ref{lemAsympLin} we show that the estimator $\widehat\theta_t$ from \eqref{eqEstimator} is asymptotically linear with influence function $x\mapsto\int\F^{-1}[\phi_\eps^{-1}(-\bull)]\ast\zeta(y)(\delta_x-\PP)(\d y)$ and thus $\widehat\theta_t$ is Gaussian regular.\par
In general, semiparametric lower bounds are constructed as the supremum of the information bounds over all regular parametric submodels. As it turns out, it suffices to apply the Cram\'er--Rao bound to the least favorable one-dimensional submodel $\PP_{g}$ of the form
\[
  f_{Y,\xi g}=f_{X,\xi g}\ast f_\eps\quad\text{ with }\quad f_{X,\xi g}:=f_X+\xi g,\quad \text{ for all }\xi\in(-\tau,\tau),
\]
with some $\tau>0$ and a perturbation $g$ satisfying
\begin{equation}\label{eqCondG}
  f_X\pm \tau g\ge0\quad\text{and}\quad\int g=0.
\end{equation}
Note that all laws $\PP_{g}$ are absolutely continuous with respect to $\PP$ assuming $\supp (f_X)=\R$. Moreover, the submodels are regular with score function $g\ast f_\eps/f_Y$, since for all $\xi\in(-\tau,\tau)\setminus\{0\}$ we have the $L^2$--differentiability
\begin{align*}
  \int\Big(\frac{f_{Y,\xi g}-f_Y-\xi g*f_\eps}{\xi f_Y}\Big)^2f_Y=0.
\end{align*}
Similarly to \citet[Chap. 25.5]{vanderVaart1998}, we define the score operator $Sg:=(g\ast f_\eps)f_Y^{-1/2}$ and thus the information operator of $f_X$ is given by $\I:=S^\star S$, where $S^\star$ denotes the adjoint of the linear operator $S$. This yields the Fisher information in direction g
\begin{equation}\label{FisherInfo}
\langle \I g,g\rangle=\scapro{Sg}{Sg}=\int\Big(\frac{g*f_\eps}{f_Y}\Big)^2f_Y
\end{equation}
and we obtain the information bound
\begin{equation}\label{eqCR}
  \mathcal I_\zeta:=\sup_g\frac{\langle g,\zeta\rangle^2}{\langle Sg,Sg\rangle},
\end{equation}
where the supremum is taken over all $g$ satisfying \eqref{eqCondG}. In the notation of \cite[Def. 3.3.2]{bicklEtAl1998}, we consider the tangent space $\dot Q:=\{(g\ast f_\eps)/f_Y|g\text{ satisfies \eqref{eqCondG}}\}$, representing the submodel $\{\PP_g\}$, and the efficient influence function of the parameter $\theta_\zeta:\dot Q\to\R, h\mapsto\scapro h\zeta$ needs to be determined.\par
Since we perturb the density additively with the restriction \eqref{eqCondG}, the quotient $|g/f_X|$ needs to be bounded and thus it is natural to assume a lower bound for the decay behavior of $f_X$. We state with some $\delta>0$ and $M\in\N$
\begin{assumption}\label{assFxPrim}\needspace{3\baselineskip}
  Let the following be satisfied
  \begin{enumerate}
    \item $f_X$ is bounded and fulfills the moment condition $\int |x|^{2+\delta}f_X(x)\d x<\infty$,
    \item $f_X\in W_1^2(\R)$ that is $f_X$ has $L^1$-Sobolev regularity two,
    \item $f_X(x)\gtrsim\langle x\rangle^{-M}$ for $x\in\R$.
  \end{enumerate}
\end{assumption}
A precise definition of the $L^1$-Sobolev space $W_1^2(\R)$ can be found in the appendix. Due to the Sobolev embedding  $W_1^2(\R)\subset H^\alpha(\R)$ with $\alpha<3/2$ (cf. \eqref{eqWpEmbed} and \eqref{eqBpEmbed}), Assumption~\ref{assFxPrim} implies the Assumption~\ref{assFx} in the previous section. The conditions on $\eps$ need to be strengthened, too.
\begin{assumption}\label{assEpsPrim}\needspace{4\baselineskip}
  We suppose
  \begin{enumerate}
    \item $\int|x|^{2+\delta} f_\eps(x)\d x<\infty$,
    \item for some $\beta\in(0,\infty)\setminus\Z$ and $M$ from above let $\phi_\eps\in C^{(\floor\beta\vee M)+1}(\R)$ satisfy for all $k=0,\dots,(\floor\beta\vee M)+1$
  \[ \mathbbm 1_{\{k=0\}}\langle u\rangle^{-\beta-k}\lesssim|\phi_\eps^{(k)}(u)|\lesssim \langle u \rangle^{-\beta-k}.\]
  \end{enumerate}
\end{assumption}
Since $M+1\ge2$, easy calculus shows that Assumption~\ref{assEps}\ref{assEpsDecay} on $\phi_{\eps}^{-1}$ follows from Assumption~\ref{assEpsPrim} on $\phi_{\eps}$. We supposed $\beta\notin\Z$ mainly to simplify our proofs. Let us first show an information bound for smooth $\zeta$.
\begin{theorem}\label{thmCRBound}
  Grant Assumptions \ref{assFxPrim} and \ref{assEpsPrim} and let $\zeta\in\mathscr S(\R)$ be a Schwartz function. For any regular estimator $T$ of $\theta_0=\langle \zeta,f_X\rangle$ with asymptotic variance $\sigma^2$ we obtain
  \begin{equation}\label{eqBound}
    \sigma^2\ge\int\big(\F^{-1}[\varphi_\eps^{-1}(-\bull)]*\zeta\big)^2f_Y -\theta_0^2.
  \end{equation}
  In particular, the supremum in \eqref{eqCR} is attained at $g^*:=g^*(\zeta):=\I ^{-1}\zeta-\scapro\zeta{f_X} f_X$, where the inverses of $S^\star$ and $\I$ are given by 
  \begin{align*}
    (S^\star)^{-1}\zeta&=(\F^{-1}[\varphi_\eps^{-1}(-\bull)]*\zeta)\sqrt{f_Y}\quad\text{and}\quad\\
    \I ^{-1}\zeta&=S^{-1}(S^{-1})^\star \zeta=\F^{-1}[\varphi_\eps^{-1}]*\big\{\big(\F^{-1}[\varphi_\eps^{-1}(-\bull)]*\zeta\big)f_Y\big\}.
  \end{align*}
\end{theorem}
Therefore, the score function corresponding to $g^*(\zeta)$ which is given by 
\[
  \F^{-1}[\varphi_\eps^{-1}(-\bull)]*\zeta -\int(\F^{-1}[\varphi_\eps^{-1}(-\bull)]*\zeta)f_Y
\]
(compare \eqref{eqSg} below) is the efficient influence function and, moreover, equals the influence function of $\widehat\theta_\zeta$. This equality shows that the estimator is efficient for smooth functionals $\theta_\zeta$. Moreover, we found already the efficient influence function in the larger tangent set of all regular submodels.\par
Unfortunately, less smooth $\zeta$ might be only in the domain of $(S^\star)^{-1}$ while $\I^{-1}\zeta$ is not in $L^2(\R)$ and thus the formal maximizer $g^*(\zeta)$ cannot be applied rigorously as the following example shows.
\begin{examples}
  Let $\eps_j$ be gamma distributed with density $\gamma_{\beta,1}$ for $\beta\in(1/4,1/2)$ and consider $\zeta(x)=e^x\mathbbm 1_{(-\infty,0]}(x)=\gamma_{1,1}(-x)$ which is contained in $Z^{\gamma_s,\gamma_c}$ for all $\gamma_s<1/2$ and $\gamma_c$ arbitrary large. We obtain 
  \[
  (S^\star)^{-1}\zeta=\gamma_{1-\beta,1}(-\bull)\sqrt{f_Y}\quad\text{and}\quad
  \I^{-1}\zeta=\F^{-1}\big[(1-iu)^\beta((1+iu)^{-1+\beta}\ast\phi)\big].
  \]
  While first term behaves nicely the Fourier transform of $\I^{-1}\zeta$ is of order $|u|^{-1+2\beta}>|u|^{-1/2}$ for $|u|\to\infty$ and thus $\I^{-1}\zeta\notin L^2(\R)$.
\end{examples}
Therefore, we choose an approximating sequence $\zeta_n\to\zeta$ with $(\zeta_n)_{n\in\N}\subset\mathscr S(\R)$. For $n\in\N$ let $g_n^*:=g^*(\zeta_n)=\I ^{-1}\zeta_n-\scapro\zeta{f_X}f_X$ be the least favorable direction in the estimation problem with respect to $\langle f_X,\zeta_n\rangle$. We obtain for every $n\in\N$
\begin{align*}
  \mathcal I_{\zeta}\ge\frac{\langle g_n^*,\zeta\rangle^2}{\langle Sg_n^*,Sg_n^*\rangle}
  = \frac{\big(\langle g_n^*,\zeta-\zeta_n\rangle+\langle g_n^*,\zeta_n\rangle\big)^2}{\langle Sg_n^*,Sg_n^*\rangle}.
\end{align*}
This inequality suggests two possibilities to understand our strategy for obtaining the efficiency bound. First, the sequence $(g_n^*)$ approximates the formal maximizer $g^*(\zeta)$ and thus plugging $g_n^*$ into the bound \eqref{eqCR} might converge to the supremum. Second, any unbiased estimator of $\theta_{\zeta_n}=\langle f_X,\zeta_n\rangle$ is at the same time a possibly biased estimator of $\theta_{\zeta}$ with bias tending to zero. Therefore, the bound for the smooth problems should converge to the nonsmooth one. The following lemma provides a sufficient condition for the convergence of the Cram\'er--Rao bounds.\par
\begin{lemma}\label{lemApprox}
  Let $\zeta$ and $(\zeta_n)$ satisfy $(S^\star)^{-1}\zeta\in L^2(\R)$ and $\zeta_n, \I ^{-1}\zeta_n\in L^2(\R)$ for all $n\in\N$. Then $\theta_{\zeta_n}\to\theta_\zeta$ and $\frac{\langle g_n^*,\zeta\rangle^2}{\langle Sg_n^*,Sg_n^*\rangle}\rightarrow\langle (S^\star)^{-1}\zeta,(S^\star)^{-1}\zeta\rangle-\scapro\zeta{f_X}^2$ hold as $n\to\infty$ if
 \[\|(S^\star)^{-1}(\zeta_n-\zeta)\|_{L^2}\to0,\quad\text{as } n\to\infty.\]
\end{lemma}
Using mapping properties on Besov spaces, we will show that the underlying Fourier multiplier $\F^{-1}[\phi_\eps^{-1}]$ and thus the inverse adjoint score operator $(S^\star) ^{-1}$ are well-defined on the set $Z^{\gamma_s,\gamma_c}$. This allows the extension of Theorem~\ref{thmCRBound} to all $\zeta\in Z^{\gamma_s,\gamma_c}$ with $\gamma_s>\beta$ and $\gamma_c>\beta+1/2$.\par
Since $\widehat\theta_t$ does not only estimate $\theta_t$ pointwise but also as a process in $\ell^\infty(\R)$, we want to generalize Theorem~\ref{thmCRBound} in this direction, too. In view of Theorem 25.48 of \citet{vanderVaart1998} the remaining ingredient is the tightness of the limiting object, which is already a necessary condition for the Donsker theorem. A regular estimator $T_n$ of $(\theta_t)_{t\in\R}$ in $\ell^\infty(\R)$ is efficient if the limiting distribution of $\sqrt n(T_n-\theta)$ is a tight zero mean Gaussian process whose covariance structure is given by the information bound for the finite dimensional distributions (cf. the convolution Theorem 5.2.1 of \citep{bicklEtAl1998}). Interestingly, the class of efficient influence functions for $t\in\R$ is not Donsker as discussed above and thus there exists no efficient estimator which is asymptotically linear in $\ell^\infty(\R)$ \cite[cf.][Thm. 18.8]{kosorok2008}.
\begin{theorem}\label{thmApprox}
  Let Assumptions \ref{assFxPrim} and \ref{assEpsPrim} be satisfied as well as $\zeta\in Z^{\gamma_s,\gamma_c}$ with $\gamma_s>\beta$ and $\gamma_c>\beta+1/2$. Then the estimator $(\widehat\theta_t)_{t\in\R}$ defined in \eqref{eqEstimator} is (uniformly) efficient.
\end{theorem}
Additionally, the proof of Theorem~\ref{thmApprox} reveals the relation between the intrinsic metric $d(s,t)^2=\E[(\mathbb G_s-\mathbb G_t)^2]$ of the limit $\mathbb G$, which is essential to show tightness, and the metric $d_{I^{-1}}(s,t)^2=\langle (S^\star)^{-1}(\zeta_t-\zeta_s),(S^\star)^{-1}(\zeta_t-\zeta_s)\rangle$ which is induced by the inverse Fisher information, namely
\[
  d_{I^{-1}}(s,t)^2=d(s,t)^2+\scapro{\zeta_t-\zeta_s}{f_X}^2
\]
(cf. equations \eqref{centering} and \eqref{eqApproxZeta} below) such that both metrics are equal up to some centering term which is another way of interpreting the efficiency of $\widehat\theta_\bull$.\par

\section{Proof of the Donsker theorem}\label{secProofs}
First, we provide an auxiliary lemma, which describes the properties of the deconvolution operator $\F^{-1}[\varphi_\eps^{-1}]$.
\begin{lemma}\label{lemPsiDO}
  Grant Assumption~\ref{assEps}.
  \begin{enumerate}
    \item \label{FourierMultiplier} For all $s\in\R, p,q\in[1,\infty]$ the deconvolution operator $\F^{-1}[\phi_\eps^{-1}(-\bull)]$ is a Fourier multiplier from $B_{p,q}^s(\R)$ to $B_{p,q}^{s-\beta}(\R)$, that is the linear map
    \[B_{p,q}^s(\R)\to B_{p,q}^{s-\beta}(\R),f\mapsto\F^{-1}[\phi_\eps^{-1}(-\bull)\F f]\]
    is  bounded.
    \item\label{lemPsiDOq} For any integer $m$ strictly larger then $\beta$ we have $\F^{-1}[(1+iu)^{-m}\varphi_\eps^{-1}]\in L^1(\R)$ and if $m>\beta+1/2$ we also have $\F^{-1}[(1+iu)^{-m}\varphi_\eps^{-1}]\in L^2(\R)$.
    \item Let $\beta^+>\beta$ and $f,g\in H^{\beta^+}(\R)$. Then
    \begin{equation}\label{eqAdjoint}
      \int\big(\F^{-1}[\varphi_\eps^{-1}]\ast f\big)g=\int\big(\F^{-1}[\varphi_\eps^{-1}(-\bull)]\ast g\big)f.
    \end{equation}
    Using the kernel $K$, this equality extends to functions $g\in L^2(\R)\cup L^\infty(\R)$ and finite Borel measures $\mu$:
    \begin{equation}\label{eqAdjointSmooth}
      \int\big(\F^{-1}[\varphi_\eps^{-1}\F K_h]\ast\mu\big)g
      =\int\big(\F^{-1}[\varphi_\eps^{-1}(-\bull)\F K_h]\ast g\big)\d \mu.
    \end{equation}
  \end{enumerate}
\end{lemma}
\begin{proof}
\hspace{1em}
\begin{enumerate}
  \item Analogously to \cite{NicklReiss2012}, we deduce from Corollary 4.11 of \cite{girardiWeis2003} that $(1+iu)^{-\beta}\phi_\eps^{-1}(-u)$ is a Fourier multiplier on $B^s_{p,q}$ by Assumption~\ref{assEps}\ref{assEpsDecay}. It remains to note that $j: B_{p,q}^s(\R)\to B_{p,q}^{s-\beta}(\R), f\mapsto \F^{-1}[(1+iu)^{\beta}\F f]$ is a linear isomorphism \citep[Thm. 2.3.8]{triebel2010}.
  \item Since the gamma density $\gamma_{1,1}$ is of bounded variation, it is contained in $B^1_{1,\infty}(\R)$ by \eqref{eqBoundVar}. Using the isomorphism $j$ from \ref{FourierMultiplier}, we deduce $\gamma_{m,1}\in B^m_{1,\infty}(\R)$ and thus by Besov embeddings \eqref{eqBqEmbed} and \eqref{eqWpEmbed}
  \[
    \F^{-1}[(1+iu)^{-m}\varphi_\eps^{-1}]\in B^{m-\beta}_{1,\infty}(\R)\subset B^0_{1,1}(\R)\subset L^1(\R).
  \]
  If $m-\beta>1/2$ we can apply the embedding $B^{m-\beta}_{1,\infty}(\R)\subset B^{m-\beta-1/2}_{2,\infty}(\R)\subset L^2(\R)$.
  \item For $f\in H^{\beta^+}(\R)$ \ref{FourierMultiplier} and the Besov embeddings \eqref{eqWpEmbed}, \eqref{eqBpEmbed} and \eqref{eqBqEmbed} yield
  \begin{align*}
    \|\F^{-1}[\varphi_\eps^{-1}]\ast f\|_{L^2}
    \lesssim\|\F^{-1}[\varphi_\eps^{-1}]\ast f\|_{B_{2,1}^0}
    \lesssim\|f\|_{B_{2,1}^\beta}
    \lesssim\|f\|_{H^{\beta^+}}<\infty.
  \end{align*}
  Therefore, it follows by Plancherel's equality
  \begin{align*}
    \int\big(\F^{-1}[\varphi_\eps^{-1}]*f\big)(x)g(x)\d x
    &=\frac{1}{2\pi}\int\varphi_\eps^{-1}(-u)\F f(-u)\F g(u)\d u\\
    &=\int\big(\F^{-1}[\varphi_\eps^{-1}(-\bull)]*g\big)(x)f(x)\d x.
  \end{align*}
  To prove the second part of the claim for $g\in L^2(\R)$, we note that by Young's inequality
  \[
    \|\F^{-1}[\phi_\eps^{-1}\F K_h]\|_{L^2}\le\|\F^{-1}[\phi_\eps^{-1}\mathbf 1_{[-1/h,1/h]}]\|_{L^2}\|K_h\|_{L^1}<\infty
  \]
  due to the support of $\F K$ and Assumption \eqref{assKDecay} on the decay of $K$. Since $\mu$ is a finite measure and $g$ is bounded, Fubini's theorem yields then
  \begin{align*}
    &\int g(x)\big(\F^{-1}[\phi_\eps^{-1}\F K_h]\ast\mu\big)(x)\d x\\
    &\quad=\int\int g(x)\F^{-1}[\phi_\eps^{-1}\F K_h](x-y)\mu(\d y)\d x\\
    &\quad=\int\big(\F^{-1}[\phi_\eps^{-1}(-\bull)\F K_h]\ast g\big)(y)\mu(\d y),
  \end{align*}
  where we have used the symmetry of the kernel. In order to apply Fubini's theorem for the case $g\in L^\infty(\R)$, too, we have to show that
  \(\|\F^{-1}[\phi_\eps^{-1}\F K_h]\|_{L^1}\) is finite.
  We replace the indicator function by a function $\chi\in C^\infty(\R)$ which equals one on $[-1/h,1/h]$ and has got compact support. We estimate
  \begin{align}\|\F^{-1}[\phi_\eps^{-1}\F K_h]\|_{L^1}\le\|\F^{-1}[\phi_\eps^{-1}\chi]\|_{L^1}\|K_h\|_{L^1}.\label{SmoothedDeconvolutionL1}
  \end{align}
  Using $\phi_\eps^{-1}\chi$ is twice continuously differentiable and has got compact support we obtain
  \begin{align*}
    \|(1+x^2)\F^{-1}[\phi_\eps^{-1}\chi](x)\|_{\infty}
    &\le \|\F^{-1}[(\Id-\D^2)\phi_\eps^{-1}\chi](x)\|_{\infty}\\
    &\le \|(\Id-\D^2)\phi_\eps^{-1}\chi\|_{L^1}<\infty,
  \end{align*}
  where we denote the identity and the differential operator by $\Id$ and $\D$, respectively. This shows that \eqref{SmoothedDeconvolutionL1} is finite.\qedhere
\end{enumerate}
\end{proof}

\subsection{Convergence of the finite dimensional distributions}
As usual, we decompose the error into a stochastic error term and a bias term:
\begin{align*}
  &\widehat\theta_t-\theta_t=\widehat\theta_t-\E[\widehat\theta_t]+\E[\widehat\theta_t]-\theta_t\\
  &\quad=\int\zeta_t(x)\F^{-1}\Big[\F K_h\frac{\phi_n-\phi}{\phi_\eps}\Big](x)\d x+\int\zeta_t(x)(K_h\ast f_X(x)-f_X(x))\d x.
\end{align*}
\subsubsection{The bias}\label{secBias}
The bias term can be estimated by the standard kernel estimator argument. Let us consider the singular and the continuous part of $\zeta$ separately. Applying Plancherel's identity and H\"older's inequality, we obtain
\begin{align*}
   &\int|\zeta^s_t(x)(K_h\ast f_X(x)-f_X(x))|\d x\\
  &=\frac1{2\pi}\int|\F\zeta^s_t(u)(\F K(hu)-1)\F f_X(-u)|\d u\\
  &\le\|\langle u\rangle^{-(\alpha+\gamma_s)}(\F K(hu)-1)\|_\infty\int\langle u\rangle^{\alpha+\gamma_s}|\F\zeta^s(u)\F f_X(u)|\d u\\
  &\le h^{\alpha+\gamma_s}\|u^{-(\alpha+\gamma_s)}(\F K(u)-1)\|_\infty\|\zeta^s\|_{H^{\gamma_s}}\|f_X\|_{H^\alpha}
\end{align*}
The term $\|u^{-(\alpha+\gamma_s)}(\F K(u)-1)\|_\infty$ is finite using the a Taylor expansion of $\F K$ around 0 with $(\F K)^{(l)}=0$ for $l=1,\dots,\lfloor \alpha+\gamma_s\rfloor$ by the order of the kernel \eqref{eqKorder}.\par
For the smooth part of $\zeta_t$ Plancherel's identity yields
\begin{align*}
  &\int|\zeta_t^c(x)(K_h\ast f_X-f_X)(x)|\d x\\
  &=\frac1{2\pi}\int|\F[\tfrac{1}{ix+1}\zeta_t^c(x)](\Id+\D)\{(\F K(hu)-1)\F f_X(-u)\}|\d u\\
  &\le\int|\F[\tfrac{1}{ix+1}\zeta_t^c(x)](\F K(hu)-1+h\F[ixK](hu))\F f_X(-u)|\d u\\
  &\quad-\int|\F[\tfrac{1}{ix+1}\zeta_t^c(x)](\F K(hu)-1)\F[ixf_X](-u)|\d u.
\end{align*}
The first term can be estimated as before and for the second term we note that $xf_X(x)\in L^2(\R)=H^0(\R)$ by Assumption~\ref{assFx}\ref{assFxMoment} such that the additional smoothness of $\frac{1}{ix+1}\zeta^c(x)$ yields the right order. Therefore, we have $|\E[\widehat\theta_t]-\theta_t|\lesssim h^{\alpha+\gamma_s}$ and thus by the choice of $h$, the bias term is of order $o(n^{-1/2})$.

\subsubsection{The stochastic error}
We notice that $\|\zeta^c-a\|_{H^{\gamma_c}}\lesssim\|\langle x\rangle^{-\tau}\|_{C^s}\|\langle x\rangle^{\tau}(\zeta^c(x)-a(x))\|_{H^{\gamma_c}}<\infty$ for any $s>\gamma_c$, where we used the pointwise multiplier property \eqref{eqPointMult} as well as the Besov embeddings \eqref{eqBqEmbed} and \eqref{eqCsEmbed}.
We have $\zeta^s\in L^2$ and by \eqref{eqWpEmbed}, \eqref{eqBpEmbed} and \eqref{eqBqEmbed}
\[\|\zeta^c\|_\infty \le \|a\|_\infty+\|\zeta^c-a\|_\infty\le \|a\|_\infty+\|\zeta^c-a\|_{H^{\gamma_c}}<\infty,\]
since $\gamma_c>1/2$.
Consequently we can apply the smoothed adjoint equality \eqref{eqAdjointSmooth} and obtain for the stochastic error term
\begin{align}
  &\int\zeta_t(x)\F^{-1}\Big[\F K_h\frac{\phi_n-\phi}{\phi_\eps}\Big](x)\d x\notag\\
  &\quad=\int\F^{-1}[\phi_\eps^{-1}(-\bull)\F K_h]\ast\zeta_t(x)(\PP_n-\PP)(\d x).\label{stochError}
\end{align}
Therefore, it suffices for the convergence of the finite dimensional distributions to bound the term
\begin{equation}\label{eqMomentCond}
  \sup_{h\in(0,1)}\int\left|\F^{-1}[\phi_\eps^{-1}(-\bull)\F K_h]\ast\zeta(x)\right|^{2+\delta}\PP(\d x),
\end{equation}
for any function $\zeta\in Z^{\gamma_s,\gamma_c}$. Then the stochastic error term converges in distribution to a normal random variable by the central limit theorem under the Lyapunov condition \citep[i.e.,][Thm.~15.43 together with Lem.~15.41]{klenke2007}. Finally, the Cram\'er-Wold device yields the convergence of the finite dimensional distributions in Theorem~\ref{thmDonsker}.\par
First, note that the moment conditions in Assumptions~\ref{assFx} and \ref{assEps} and the estimate
\[
  |x|^p f_Y(x)\le \int|x-y+y|^pf_X(x-y)f_\eps(y)\d y\lesssim (|y|^pf_X)*f_\eps+f_X*(|y|^pf_\eps),
\]
for $x\in\R$, $p\ge1$, yield finite $(2+\delta)$th moments for $\PP$ since
\begin{equation}\label{momentfy}
  \int |x|^{2+\delta}f_Y(x)\d x
  \lesssim \| |x|^{2+\delta}f_X\|_{L^1}\|f_\eps\|_{L^{1}}+\|f_X\|_{L^1}\| |x|^{2+\delta}f_\eps\|_{L^{1}}<\infty.
\end{equation}
To estimate \eqref{eqMomentCond}, we rewrite
\begin{align}
  \F^{-1}[\varphi_\eps^{-1}(-\bull)]\ast\zeta^c(x)
  &=\F^{-1}\big[\varphi_\eps^{-1}(-u)(\Id+\D)\F[\tfrac{1}{iy+1}\zeta^c(y)](u)\big](x)\notag\\
  &=\F^{-1}\big[\varphi_\eps^{-1}(-u)\F[\tfrac{1}{iy+1}\zeta^c(y)](u)\big](x) \notag\\ &\quad+\F^{-1}\big[\varphi_\eps^{-1}(-u)\big(\F[\tfrac{1}{iy+1}\zeta^c(y)]\big)'(u)\big](x)\label{eqDerivRepr}\\
  &=(1+ix)\F^{-1}[\varphi_\eps^{-1}(-u)\F[\tfrac{1}{iy+1}\zeta^c(y)](u)](x) \notag\\ &\quad+\F^{-1}[(\varphi_\eps^{-1})'(-u)\F[\tfrac{1}{iy+1}\zeta^c(y)](u)](x),\notag
\end{align}
owing to the product rule for differentiation. Hence,
\begin{align}
  \F^{-1}[\varphi_\eps^{-1}(-\bull)]\ast\zeta(x)
  =&\F^{-1}[\varphi_\eps^{-1}(-u)\F \zeta^s(u)](x)\notag\\
  &\quad+(1+ix)\F^{-1}[\varphi_\eps^{-1}(-u)\F[\tfrac{1}{iy+1}\zeta^c(y)](u)](x)\notag\\
  &\quad+\F^{-1}[(\varphi_\eps^{-1})'(-u)\F[\tfrac{1}{iy+1}\zeta^c(y)](u)](x).\label{eqDecomp}
\end{align}
While $\F^{-1}[\varphi_\eps^{-1}(-\bull)]\ast\zeta$ may exist only in distributional sense in general, it is defined rigorously through the right-hand side of the above display for $\zeta\in Z^{\gamma_s,\gamma_c}$. Considering $\zeta\ast K_h$ instead of $\zeta$, we estimate separately all three terms in the following.\par
The continuity and linearity of the Fourier multiplier $\F^{-1}[\phi_\eps^{-1}(-\bull)]$, which was shown in Lemma~\ref{lemPsiDO}\ref{FourierMultiplier}, yield for the first term in \eqref{eqDecomp}
\begin{align*}
  \|\F^{-1}[\varphi_\eps^{-1}(-u)\F\zeta^s(u)\F K_h(u)]\|_{H^{\delta}}
  &=\big\|\F^{-1}\big[\varphi_\eps^{-1}(-\bull)\F[\zeta^s\ast K_h]\big]\big\|_{B^{\delta}_{2,2}}\\
  &\lesssim\|\zeta^s*K_h\|_{B^{\beta+\delta}_{2,2}}
  \lesssim\|\zeta^s\|_{H^{\beta+\delta}},
\end{align*}
where the last inequality holds by $\|\F K_h\|_\infty\le\|K\|_{L^1}$. Using the boundedness of $f_Y$ and the continuous Sobolev embedding $H^{\delta/4}(\R)\subset L^{2+\delta}(\R)$ by \eqref{eqWpEmbed}, \eqref{eqBqEmbed} and \eqref{eqBpEmbed}, we obtain
\begin{align}
  &\|\F^{-1}[\varphi_\eps^{-1}(-u)\F\zeta^s(u)\F K_h(u)]\|_{L^{2+\delta}(\PP)}\notag\\
  &\quad\lesssim\|\F^{-1}[\varphi_\eps^{-1}(-u)\F\zeta^s(u)\F K_h(u)]\|_{L^{2+\delta}}\notag\\
  &\quad\lesssim\|\F^{-1}[\varphi_\eps^{-1}(-u)\F\zeta^s(u)\F K_h(u)]\|_{H^{\delta}}\notag\\
  &\quad\lesssim\|\zeta^s\|_{H^{\beta+\delta}}\label{eq2+momSing}
\end{align}
To estimate the second term in \eqref{eqDecomp}, we use the Cauchy--Schwarz inequality and Assumption~\ref{assEps}\ref{assEpsDecay}:
\begin{align*}
  &\|\F^{-1}[\varphi_\eps^{-1}(-u)\F[\tfrac{1}{ix+1}\zeta^c(x)](u)\F K_h(u)]\|_\infty\\
  &\qquad\le\|\varphi_\eps^{-1}(-u)\F[\tfrac{1}{ix+1}\zeta^c]\F K_h(u)\|_{L^1}\\
  &\qquad\lesssim\|\langle u\rangle^{-1/2-\beta-\delta}\varphi_\eps^{-1}(-u)\|_{L^2}\|\langle u\rangle^{1/2+\beta+\delta}\F[\tfrac{1}{ix+1}\zeta^c(x)]\|_{L^2}\\
  &\qquad\lesssim\|\tfrac{1}{ix+1}\zeta^c(x)\|_{H^{1/2+\beta+\delta}}.
\end{align*}
Thus $\int (1+x^2)^{(2+\delta)/2}f_Y(x)\d x<\infty$ from \eqref{momentfy} yields
\begin{align}
  &\|(1+ix)\F^{-1}[\varphi_\eps^{-1}(-u)\F[\tfrac{1}{iy+1}\zeta^c(y)](u)\F K_h(u)](x)\|_{L^{2+\delta}(\PP)}\notag\\
  &\qquad\lesssim\|\tfrac{1}{ix+1}\zeta^c(x)\|_{H^{1/2+\beta+\delta}}.\label{eq2+momCont1}
\end{align}
The last term in the decomposition \eqref{eqDecomp} can be estimated similarly using the Cauchy--Schwarz inequality and Assumption~\ref{assEps}\ref{assEpsDecay} for $(\phi^{-1})'$
\begin{align}
  &\|\F^{-1}[(\varphi_\eps^{-1})'(-u)\F[\tfrac{1}{ix+1}\zeta^c(x)](u)\F K_h(u)]\|_{L^{2+\delta}(\PP)}\notag\\
  &\qquad\lesssim\|(\varphi_\eps^{-1})'(-u)\F[\tfrac{1}{ix+1}\zeta^c(x)](u)\|_{L^1}\notag\\
  &\qquad\le\|\langle u\rangle^{1/2-\beta-\delta}(\varphi_\eps^{-1})'\|_{L^2}\|\langle u\rangle^{-1/2+\beta+\delta}\F^{-1}[\tfrac{1}{ix+1}\zeta^c(x)](u)\|_{L^2}\notag\\
  &\qquad\lesssim\|\tfrac{1}{ix+1}\zeta^c(x)\|_{H^{-1/2+\beta+\delta}}.\label{eq2+momCont2}
\end{align}
Combining \eqref{eq2+momSing}, \eqref{eq2+momCont1} and \eqref{eq2+momCont2}, we obtain
\begin{equation}\label{eq2+bound}
  \sup_{h\in(0,1)}\|\F^{-1}[\phi_\eps^{-1}(-\bull)\F K_h]\ast\zeta(x)\|_{L^{2+\delta}(\PP)}
  \lesssim\|\zeta\|_{Z^{\beta+\delta,1/2+\beta+\delta}},
\end{equation}
which is finite for $\delta$ small enough satisfying $\beta+\delta\le\gamma_s$ and $1/2+\beta+\delta\le\gamma_c$.
Since $\F K_h$ converges pointwise to one and $|\F^{-1}[\phi_\eps^{-1}(-\bull)\F K_h]\ast\zeta(x)|^2$ is uniformly integrable by the bound of the $2+\delta$ moments, the variance converges to
\[
\int\left|\F^{-1}[\phi_\eps^{-1}(-\bull)]\ast\zeta(x)\right|^2\PP(\d x).
\]

\subsection{Tightness}\label{secTightness}
Motivated by the representation \eqref{stochError} of the stochastic error,
we introduce the empirical process
\begin{align}\label{defEmpiricalProcess}
  \nu_n(t):=\sqrt{n}\int\F^{-1}[\phi_\eps^{-1}(-\bull)\F K_h]\ast\zeta_{t}(x)(\PP_n-\PP)(\d x),\quad t\in\R.
\end{align}
In order to show tightness of the empirical process, we first show some properties of the class of translations $\mathcal H:=\{\zeta_t|t\in\R\}$ for $\zeta\in Z^{\gamma_s,\gamma_c}$.
\begin{lemma}\label{lemTranslation}
  For $\zeta\in Z^{\gamma_s,\gamma_c}$ the following is satisfied:
  \begin{enumerate}
    \item\label{uniformZnorm} The decomposition $\zeta_t=\zeta^c_t+\zeta^s_t$ satisfies the conditions in the definition of $Z^{\gamma_s,\gamma_c}$ with $a_t$. We have $\sup_{t\in\R} \|\zeta_t\|_{Z^{\gamma_s,\gamma_c}}<\infty$.
    \item For any $\eta\in(0,\gamma_s)$ there is a $\tau>0$ such that $\|\zeta_t-\zeta_s\|_{Z^{\gamma_s-\eta,\gamma_c-\eta}} \lesssim |t-s|^\tau$ holds for all $s,t\in \R$ with $|t-s|\le1$
  \end{enumerate}
\end{lemma}
\begin{proof}\hspace{1em}
  \begin{enumerate}
    \item Since $\|\zeta^s_t\|_{H^{\gamma_s}}^2=\int\langle u\rangle^{2\gamma_s}|e^{itu}\F\zeta^s(u)|^2\d u=\|\zeta^s\|_{H^{\gamma_s}}^2$, both claims hold for the singular part.
        Applying the pointwise multiplier property of Besov spaces
        \eqref{eqPointMult}  as well as the Besov embeddings \eqref{eqBqEmbed} and \eqref{eqCsEmbed}, we obtain for some $M>\gamma_c$ and $a\in C^\infty(\R)$ as in definition~\eqref{eqZ}
    \begin{align*}
          \|\langle x \rangle ^\tau\big(\zeta^c_t(x)-a_t(x)\big)\|_{H^{\gamma_c}} &\lesssim \|\tfrac{\langle x \rangle ^\tau}{\langle x-t \rangle ^\tau}\|_{C^M} \|\langle x-t \rangle ^\tau\big(\zeta^c_t(x)-a_t(x)\big)\|_{H^{\gamma_c}}\\
          &= \|\tfrac{\langle x \rangle ^\tau}{\langle x-t \rangle ^\tau}\|_{C^M} \|\langle x \rangle ^\tau\big(\zeta^c(x)-a(x)\big)\|_{H^{\gamma_c}},
    \end{align*}
    which is finite for all $t\in\R$ since $\langle x \rangle ^\tau\langle x-t \rangle ^{-\tau}\in C^M(\R)$. For the second claim we estimate similarly
    \begin{align*}
      \sup_{t\in\R}\|\tfrac{1}{ix+1}\zeta_t^c(x)\|_{H^{\gamma_c}} &\lesssim\sup_{t\in\R}\|\tfrac{a_t(x)}{ix+1}\|_{H^{\gamma_c}} +\|\tfrac{1}{ix+1}\|_{C^M}\sup_{t\in\R}\|\zeta_t^c-a_t\|_{H^{\gamma_c}}\\
      &\lesssim \|\tfrac{1}{ix+1}\|_{H^{\gamma_c}} \|a\|_{C^M} +\|\tfrac{1}{ix+1}\|_{C^M}\|\zeta^c-a\|_{H^{\gamma_c}}<\infty.
    \end{align*}
    \item For the singular part note that
    \begin{align*}
      &\|\zeta^s_t-\zeta^s_s\|_{H^{\gamma_s-\eta}}\\
      &\le\|\langle u\rangle^{\gamma_s}\F \zeta^s(u)\|_{L^2}\|\langle u\rangle^{-\eta}(1-e^{i(t-s)u})\|_\infty\\
      &\lesssim \|\langle u\rangle^{-\eta}\|_{L^\infty(\R\backslash(-|t-s|^{-1/2},|t-s|^{-1/2}))}\\
      &\qquad\vee \|(1-e^{i(t-s)u})\|_{L^\infty((-|t-s|^{-1/2},|t-s|^{-1/2}))}\\
      &\lesssim |t-s|^{\eta/2}\vee|t-s|^{1/2}.
    \end{align*}
    For $\zeta^c$ we have
    \begin{align*}  \left\|\tfrac{1}{ix+1}(\zeta^c_t(x)-\zeta^c_s(x))\right\|_{H^{\gamma_c-\eta}}
    &\lesssim \left\|\tfrac{1}{ix+1}\zeta^c_t(x)-\big(\tfrac{1}{ix+1}\zeta^c_t(x)\big)\ast\delta_{s-t} \right\|_{H^{\gamma_c-\eta}} \\ &\quad+\left\|\tfrac{1}{i(x-s+t)+1}\zeta^c_s(x)-\tfrac{1}{ix+1}\zeta^c_s(x)\right\|_{H^{\gamma_c-\eta}}.
	\end{align*}
    The first term can be treated analogously to $\zeta^s$. Using some integer $M\in\N$ strictly larger than $\gamma_c$, the second term can be estimated by
    \begin{align*}
    &\left\|\tfrac{1}{i(x-s+t)+1}\zeta^c_s(x)-\tfrac{1}{ix+1}\zeta^c_s(x)\right\|_{H^{\gamma_c-\eta}}\\
    &\quad\lesssim |t-s| \left\| \tfrac{1}{i(x-s+t)+1}\tfrac{1}{ix+1}\zeta^c_s(x)\right\|_{H^{\gamma_c-\eta}}\\
    &\quad\lesssim|t-s| \left\| \tfrac{1}{i(x-s+t)+1}\right\|_{C^{M}}\left\|\tfrac{1}{ix+1}\zeta^c_s(x)\right\|_{H^{\gamma_c-\eta}}\\
    &\quad\lesssim |t-s|,
    \end{align*}
    where we used again pointwise multiplier \eqref{eqPointMult}, embedding  properties of Besov spaces \eqref{eqBqEmbed} and \eqref{eqCsEmbed} as well as~\ref{uniformZnorm}.\qedhere
  \end{enumerate}
\end{proof}

\subsubsection{Pregaussian limit process}\label{secPregaussian}

Let $\mathbb G$ be the stochastic process from Theorem~\ref{thmDonsker}. It induces the intrinsic covariance metric $d(s,t):=\E[(\mathbb G_s-\mathbb G_t)^2]^{1/2}$.

\begin{theorem}\label{thmPregauss}
  There exists a version of $\mathbb G$ with uniformly $d$-continuous sample paths almost surely and with $\sup_{t\in\R} |\mathbb G_t|<\infty$ almost surely.
\end{theorem}

The proof of the theorem shows in addition that~$\R$ is totally bounded with respect to~$d$. The boundedness of the sample paths follows from the totally bounded index set and the uniform continuity.
Further we conclude that $\G$ defined in~\eqref{g} is $\PP$--pregaussian by \citet[p.~89]{vanderVaartWellner1996}. Thus~$\mathbb G$ is a tight Borel random variable in $\ell^\infty(\R)$ and the law of~$\mathbb G$ is uniquely defined through the covariance structure and the sample path properties in the theorem \citep[Lem.~1.5.3]{vanderVaartWellner1996}.


\begin{proof}
To show that the class is pregaussian, it suffices to verify polynomial covering numbers. To that end we deduce that
\begin{equation}\label{centering}
d(s,t)=\big(\|g_t-g_s\|_{L^2(\PP)}^2-\scapro{\zeta_t-\zeta_s}{f_X}^2\big)^{1/2}\le \|g_t-g_s\|_{L^2(\PP)}
\end{equation}
decreases polynomial for $|t-s|\to0$, for $\max(s,t)\to\infty$ and for $\min(s,t)\to\infty$.
Using the same estimates which show the moment bound \eqref{eq2+bound} but replacing $\F K_h=1$, we obtain
\begin{equation}\label{eqMomentBound}
  \|\F[\varphi_\eps^{-1}(-\bull)]*\zeta\|_{L^2(\PP)}
  \lesssim\|\zeta\|_{Z^{\beta+\delta,1/2+\beta+\delta}}
\end{equation}
and thus by choosing $\delta$ and $\eta$ small enough Lemma~\ref{lemTranslation} yields
$d(s,t)\lesssim\|\zeta_t-\zeta_s\|_{Z^{\beta+\delta,1/2+\beta+\delta}}\lesssim|t-s|^\tau$.
We now turn to the estimation of the tails. We will only consider the case $s,t\ge N$ since the case $s,t\le N$ can be treated in the same way. Without loss of generality, let $s<t$.

For the smooth component of $\zeta$ we have to show that
$\big\|\tfrac{1}{ix+1}(\zeta^c_t(x)-\zeta^c_s(x))\big\|_{H^{\gamma_c}}$ with $t,s\ge N$ decays polynomially in $N$. It suffices to prove $\big\|\tfrac{1}{ix+1}(\zeta^c_t-a_t)(x)\big\|_{H^{\gamma_c}}$ and $\big\|\tfrac{1}{ix+1}(a_t-a_s)(x)\big\|_{H^{\gamma_c}}$ with $a\in C^\infty(\R)$ from definition \eqref{eqZ} of $Z^{\gamma_s,\gamma_c}$ both decay polynomially in $N$. Let $M>\gamma_c$ and $\psi\in C^M(\R)$ with $\psi(x)=1$ for $x\in\R\setminus[-\frac{1}{2},\frac{1}{2}]$ and $\psi(x)=0$ for $x\in[-\frac{1}{4},\frac{1}{4}]$. The pointwise multiplier property \eqref{eqPointMult} yields
\begin{align*}
  &\big\|\tfrac{1}{ix+1}(\zeta^c_t-a_t)(x)\big\|_{H^{\gamma_c}}\\
  &\quad=\big\|\big(\psi(x/N)+(1-\psi(x/N))\big)\tfrac{1}{ix+it+1}(\zeta^c-a)(x)\big\|_{H^{\gamma_c}}\\
  &\quad\lesssim \|\tfrac{1}{ix+it+1}\|_{C^M}\|\psi(x/N)(\zeta^c-a)(x)\|_{H^{\gamma_c}}+\|\tfrac{1-\psi(x/N)}{ix+it+1}\|_{C^M}\|\zeta^c-a\|_{H^{\gamma_c}}\\
  &\quad\lesssim \|\langle x\rangle^{-\tau}\psi(x/N)\|_{C^M}\|\langle x\rangle^\tau(\zeta^c-a)(x)\|_{H^{\gamma_c}}+N^{-1}\|\zeta^c-a\|_{H^{\gamma_c}}\\
  &\quad\lesssim N^{-(\tau\wedge 1)}
\end{align*}
and for $N$ large enough such that $\supp (a')\subset [-N/2,N/2]$ we obtain
\begin{align*}
  &\big\|\tfrac{1}{ix+1}(a_t-a_s)(x)\big\|_{H^{\gamma_c}}\\
  &\quad=\big\|\tfrac{\psi(x/N)}{ix+1}(a_t-a_s)(x)\big\|_{H^{\gamma_c}}
  \lesssim \big\|\tfrac{\psi(x/N)}{ix+1}\big\|_{H^{\gamma_c}}\big\|(a_t-a_s)(x)\big\|_{C^M}\\
  &\quad\lesssim \big\|(ix+1)^{-3/4}\big\|_{H^{\gamma_c}} \big\|\psi(x/N)(ix+1)^{-1/4}\big\|_{C^M}
  \lesssim N^{-1/4}
  .
\end{align*}
To bound the singular part it suffices to show that
\begin{align*}
  \left\|\F^{-1}[\phi_\epsilon^{-1}(-\bull)]\ast\zeta^s_t\right\|_{L^2(\PP)},\qquad t\ge N,
\end{align*}
decays polynomially in $N$. To this end, we split the integral domain into
\begin{align}
  \left\|\F^{-1}[\phi_\epsilon^{-1}(-\bull)]\ast\zeta^s_t\right\|_{L^2(\PP)}^2
  =& \int_{-\infty}^{-N/2}|\F^{-1}[\phi_\epsilon^{-1}(-\bull)\F\zeta^s](x)|^2 f_Y(x+t)\d x\notag\\ &\quad+\int_{-N/2}^{\infty}|\F^{-1}[\phi_\epsilon^{-1}(-\bull)\F\zeta^s](x)|^2 f_Y(x+t)\d x.\label{singularTails}
\end{align}
To estimate the first term, we use the following auxiliary calculations
\begin{align*}
  &ix\F^{-1}[\phi_\eps^{-1}(-\bull)\F\zeta^s](x)\\
  &\quad=-\F^{-1}[(\phi_\eps^{-1})'(-\bull)\F\zeta^s](x) +\F^{-1}[\phi_\eps^{-1}(-\bull)\F[iy\zeta^s(y)]](x)
\end{align*}
and with an integer $M\in\N$ strictly larger than $\gamma_s$ and a function $\chi\in C^{M}(\R)$ which is equal to one on $\supp(\zeta^s)$ and has compact support
\begin{align*}
  \|y\zeta^s(y)\|_{H^{\gamma_s}}=\|y\chi(y)\zeta^s(y)\|_{H^{\gamma_s}}&\lesssim \|y\chi(y)\|_{B_{\infty,2}^{\gamma_s}}\|\zeta^s(y)\|_{H^{\gamma_s}}\\
  &\lesssim \|y\chi(y)\|_{C^{M}}<\infty,
\end{align*}
where we used the pointwise multiplier property \eqref{eqPointMult} of Besov spaces as well as the Besov embeddings \eqref{eqBqEmbed} and \eqref{eqCsEmbed}.
Thus $ix\F^{-1}[\phi_\eps^{-1}(-\bull)\F\zeta^s](x)\in L^2(\R)$. Applying this and the boundedness of $f_Y$ to the first term in \eqref{singularTails} yields
\begin{align*}
  &\int_{-\infty}^{-N/2}|\F^{-1}[\phi_\epsilon^{-1}(-\bull)\F\zeta^s](x)|^2 f_Y(x+t)\d x\\
  &\quad\lesssim \int_{-\infty}^{-N/2}|\F^{-1}[\phi_\epsilon^{-1}(-\bull)\F\zeta^s](x)|^2\d x\\
  &\quad
  \le 4N^{-2} \int_{-\infty}^{-N/2}|x\F^{-1}[\phi_\epsilon^{-1}(-\bull)\F\zeta^s](x)|^2 \d x\lesssim N^{-2}.
\end{align*}
Using H\"olders's inequality and the boundedness of $f_Y$, we estimate the second term in \eqref{singularTails} by
\begin{align*}
  &\|\F^{-1}[\phi_\epsilon^{-1}(-\bull)\F\zeta^s](x)\|_{L^{2+\delta}}^2\left(\int_{-N/2}^\infty |f_Y(x+t)|^{(2+\delta)/\delta}\d x\right)^{\delta/(2+\delta)}\\
  &\lesssim\|\F^{-1}[\phi_\epsilon^{-1}(-\bull)\F\zeta^s](x)\|_{L^{2+\delta}}^2\left(\int_{N/2}^\infty f_Y(x)\d x\right)^{\delta/(2+\delta)}.
\end{align*}
While the first factor is finite according to our bound \eqref{eq2+momSing}, which also holds when $\F K_h$ is omitted, the second one is of order $N^{-\delta}$ due to the finite $(2+\delta)$th moment of $\PP$. Therefore, the second term in \eqref{singularTails} decays polynomially.
\end{proof}

\subsubsection{Uniform central limit theorem}\label{secCLT}

We recall the definition of the empirical process $\nu_n$ in \eqref{defEmpiricalProcess}.

\begin{theorem}\label{thmTightness}
  Grant Assumptions~\ref{assFx} and \ref{assEps}.
  Let
  \[(\nu_n(t_1),\dots,\nu_n(t_k))\overset{\mathcal L}{\longrightarrow}(\mathbb G_{t_1},\dots,\mathbb G_{t_k})\]
  for all $t_1,\dots,t_k\in\R$ and for all $k\in\N$. If either $\gamma_s\le\beta+1/2$ and $h_n^\rho n^{1/4}\to\infty$ as $n\to\infty$ for some $\rho>\beta-\gamma_s+1/2$
  or if $\gamma_s>\beta+1/2$,
  then
  \[\nu_n\overset{\mathcal L}{\longrightarrow}\mathbb G \quad\text{in }\ell^\infty(\R).\]
\end{theorem}

\begin{proof}
  We split the empirical process $\nu_n$ into three parts
  \[\nu_n=\sqrt{n}\int(T_1(x)+T_2(x)+T_3(x))(\PP_n-\PP)(\d x),\]
  where $T_1$, $T_2$ and $T_3$ correspond to the three terms in decomposition~\eqref{eqDecomp} and are given by \eqref{T1}, \eqref{T2} and \eqref{T3} below.
  For the first term
  \begin{equation}\label{T1}
  T_1(x)=\F^{-1}[\varphi_\eps^{-1}(-u)\F \zeta_t^s(u)\F K_h(u)](x)
  \end{equation}
  we distinguish the two cases $\gamma_s>\beta+1/2$ and $\gamma_s\le\beta+1/2$. In the first case we will show that $T_1$ varies in a fixed Donsker class. In the second case the process indexed by $T_1$ is critical, this is where smoothed empirical processes and the condition on the bandwidth are needed. Tightness of $T_1$ in this case will be shown in Section~\ref{secCritical}. We will further show that the second term~$T_2$ and the third term~$T_3$ are both varying in fixed Donsker classes for all $\gamma_s>\beta$. In particular the three processes indexed by $T_1$, $T_2$ and $T_3$, respectively, are tight. Applying the equicontinuity characterization of tightness \citep[Thm.~1.5.7]{vanderVaartWellner1996} with the maximum of the semimetrics yields that $\nu_n$ is tight. Since we have assumed convergence of the finite dimensional distribution the convergence of $\nu_n$ in distribution follows \citep[Thm.~1.5.4]{vanderVaartWellner1996}.

  Here we consider only the first case $\gamma_s>\beta+1/2$. We recall that $\zeta^s_t$ is contained in $H^{\gamma_s}(\R)$. By the Fourier multiplier property of the deconvolution operator in Lemma~\ref{lemPsiDO}\ref{FourierMultiplier} and by $\sup_{h>0,u}|\F K_h(u)| \le \|K\|_{L^1}<\infty$ the functions $T_1$ are contained in a bounded set of $H^{1/2+\eta}(\R)$ for some $\eta>0$ small enough. We apply \citep[Prop.~1]{NicklPoetscher2007} with $p=q=2$ and $s=1/2+\eta$ and conclude that $T_1$ varies in a universal Donsker class.

  The second term is of the form
  \begin{equation}\label{T2}T_2(x)=(1+ix)\F^{-1}[\phi_\eps^{-1}(-u)\F[\tfrac{1}{iy+1}\zeta_t^c(y)](u)\F K_h(u)](x).\end{equation}
  By Assumption~\ref{assEps}\ref{assEpsDecay} we have $\phi_\eps^{-1}(u)\lesssim \langle u\rangle^{\beta}$.
  For some $\eta>0$ sufficiently small, the functions $\tfrac{1}{iy+1}\zeta_t^c(y)$, $t\in\R$, are contained in a bounded set of $H^{\beta+\eta+1/2}(\R)$ by Lemma~\ref{lemTranslation}. We obtain that the functions $T_2(x)/(1+ix)$ are contained in a bounded subset of $H^{1/2+\eta}(\R)$. Corollary~5 in \citep{NicklPoetscher2007} yields with $p=q=2$, $\beta=-1$, $s=1/2+\eta$ and $\gamma=\eta$ that $T_2$ is contained in a fixed  $\PP$-Donsker class.

  Similarly, we treat the third term
  \begin{equation}\label{T3}T_3(x)=\F^{-1}[(\phi_\eps^{-1})'(-u)\F[\tfrac{1}{iy+1}\zeta_t^c(y)](u)\F K_h(u)](x).\end{equation}
  By Assumption~\ref{assEps}\ref{assEpsDecay} we have $(\phi_\eps^{-1})'\lesssim\langle u\rangle^{\beta-1}$. As above we conclude that the functions $T_3$ are contained in a bounded set of $H^{\eta+3/2}(\R)$. By \citep[Prop.~1]{NicklPoetscher2007} with $p=q=2$ and $s=\eta+3/2$ the term $T_3$ varies in a universal Donsker class.
\end{proof}

\subsubsection{The critical term}\label{secCritical}

In this section, we treat the first term $T_1$ in the case $\gamma_s\le\beta+1/2$.
We define
\begin{equation}\label{q_t}
q_t:=\F^{-1}[\phi_\eps^{-1}(-u)\F\zeta_t^s(u)].
\end{equation}
For simplicity in point (e) below it will be convenient to work with functions $K_h$ of bounded support.
Thus we fix $\xi>0$ and define the truncated kernel
\[K_h^{(0)}:=K_h\mathbf{1}_{[-\xi,\xi]}.\]
By the Assumption \eqref{assKDecay} on the decay of $K$ we have $\sup_{h>0}\|K_h-K_h^{(0)}\|_{BV}<\infty$. We conclude $\F(K_h-K_h^{(0)})(u) \lesssim (1+|u|)^{-1}$ with a constant independent of $h>0$. By Assumption~\ref{assEps}\ref{assEpsDecay} we have $|\varphi_\eps^{-1}(u)|\lesssim(1+|u|)^{\beta}$. The functions $\zeta_t^s(u)$, $t\in\R$, are contained in a bounded set of $H^{\gamma_s}(\R)$. Consequently $T_1$ with $K_h-K^{(0)}_h$ instead of $K_h$ is contained in a bounded set of $H^{\gamma_s-\beta+1}(\R)$. With the same argument as used for $T_3$ we see that this term is contained in a universal Donsker class because $\gamma_s-\beta+1>1$ by assumption. So it remains to consider $T_1$ with the truncated kernel $K_h^{(0)}$.

In order to show tightness of the process indexed by $T_1$ with the truncated kernel $K_h^{(0)}$ we check the assumptions of Theorem~3 by \citet{GineNickl2008} in the version of~\citet[Thm.~12]{NicklReiss2012} for the class $Q=\{q_t|t\in\R\}$ and for $\mu_n(\d x):=K_{h_n}^{(0)}(x)\d x$, where $q_t(x)$ was defined in~\eqref{q_t}.
By Section~\ref{secPregaussian} the class $\G$ is $\PP$-pregaussian. From the proof also follows that $Q$ is $\PP$-pregaussian since this is just the case $\zeta^c=0$.

We write
\[Q'_\tau:=\big\{r-q\big|r,q\in Q, \|r-q\|_{L^2(\PP)}\le\tau\big\}.\]
Let $\rho>\beta-\gamma_s+1/2\ge0$ be such that $h_n^\rho n^{1/4}\to\infty$. We fix some $\rho'\in(\beta-\gamma_s+1/2,\rho \wedge 1)$ and obtain
 $h_n^{\rho'}\log(n)^{-1/2}n^{1/4}\to\infty$. We need to verify the following conditions.

\begin{enumerate}[(a)]
\item  We will show that the functions in $\tilde Q_n:=\{q_t\ast\mu_n|t\in\R\}$ are bounded by $M_n:=C h_n^{-\rho'}$ for some constant $C>0$. Since $q_t$ is only a translation of $q_0$ it suffices to consider $q_0$.
    By the definition of $Z^{\gamma_s,\gamma_c}$ in \eqref{eqZ}, by Lemma~\ref{lemPsiDO}\ref{FourierMultiplier} and by the Besov embedding \eqref{eqBqEmbed}
    \begin{align*}
     q_0&=\F^{-1}[\phi_\eps^{-1}(-u)\F\zeta^s(u)]\in B^{\gamma_s-\beta}_{2,2}(\R)\subset B^{1/2-\rho'}_{2,\infty}(\R).
    \end{align*}
    By our assumptions on the kernel~\eqref{assKDecay} it follows that $K'$ is integrable and thus that $K$ is of bounded variation.
Next we apply continuous embeddings for Besov spaces \eqref{eqWpEmbed} and \eqref{eqBpEmbed}, \eqref{eqConvolutionMult} as well as the estimate for $\|K_{h_n}\|_{B^{\rho'}_{1,1}}$ in \citet[p. 384]{GineNickl2008}, which also applies to truncated kernels, and obtain
		\begin{align}
		\|q_0\ast K_{h_n}^{(0)}\|_\infty & \lesssim  \|q_0\ast K_{h_n}^{(0)}\|_{B^0_{\infty,1}}
		\lesssim  \|q_0\ast K_{h_n}^{(0)}\|_{B^{1/2}_{2,1}}
		\lesssim \|K_{h_n}^{(0)}\|_{B^{\rho'}_{1,1}}
		\lesssim h_n^{-\rho'}.\label{envelopes}
		\end{align}

\item For $r\in Q'_\tau$ holds $\|r\ast K_h^{(0)}\|_{L^2(\PP)} \le \|r\ast K_h^{(0)}-r\|_{L^2(\PP)} +\tau$. Thus it suffices to show that $\| q\ast K_h^{(0)}-q\|_{L^2(\PP)} \to 0$ uniformly over $q\in Q$. We estimate
    \begin{align*}
      \|q_t\ast K_h^{(0)}-q_t\|_{L^2(\PP)}
      &\lesssim \|\phi_\epsilon^{-1}(-\bull)\F \zeta^s(\F K_h^{(0)}-1)\|_{L^2}.
    \end{align*}
    $\phi_\epsilon^{-1}(-\bull)\F \zeta^s$ is an $L^2$--function and
    $\F K_h^{(0)}$ is uniformly bounded and converges to one as $h\to0$. By dominated convergence the integral converges to zero.

\item The estimates in (a) can be used to see that the classes $\tilde Q_n$ have polynomial $L^2(\QQ)$--covering numbers, uniformly in all probability measures $\QQ$ and uniformly in $n$.
	The function $q_0\ast K_{h_n}^{(0)}$ is the convolution of two $L^2$-functions and thus continuous.
	The estimate \eqref{envelopes} and embedding \eqref{eqPVar} yield that $q_0\ast K_{h_n}^{(0)}$ is of finite 2--variation. We argue as in Lemma~1 by \citet{GineNickl2009}. As function of bounded 2--variation $q_0\ast K_{h_n}^{(0)}$ can be written as a composition $g_n\circ f_n$ of a nondecreasing function $f_n$ and a function $g_n$, which satisfies a Hölder condition $|g_n(u)-g_n(v)|\le |u-v|^{1/2}$, see, for example, \cite[p. 1971]{dudley1992}. More precisely, we can take $f_n(x)$ to be the 2--variation of $q_0\ast K_{h_n}^{(0)}$ up to $x$ and the envelopes of $f_n$ to be multiples of $M_n^2=C^2h_n^{-2\rho'}$.
    The set $F_n$ of all translates of the nondecreasing function $f_n$ has VC--index 2 and thus polynomial $L^1(\QQ)$--covering numbers \cite[Thm.~5.1.15]{delaPenaGine1999}. Since each $\epsilon^2$--covering of translates of $f_n$ for $L^1(\QQ)$ induces an $\epsilon$--covering of translates of $g_n\circ f_n$ for $L^2(\QQ)$ we can estimate the covering numbers by
    \begin{align*}
    N(\tilde Q_n,L^2(\QQ),\epsilon)\le N(F_n,L^1(\QQ),\epsilon^2)\lesssim (M_n/\epsilon)^4,
    \end{align*}
    with constants independent of $n$ and $\QQ$.
    The conditions for inequality (22) by \citet{GineNickl2008} are fulfilled, where the envelopes are $M_n=C h_n^{-\rho'}$ and $H_n(\eta)=H(\eta)=C_1\log(\eta)+C_0$ with $C_0,C_1>0$. Consequently
    \begin{align*}
    \E^\ast\left\|\frac1{\sqrt{n}}\sum_{j=1}^n\eps_j f(X_j)\right\|_{(\tilde Q_n)_{n^{-1/4}}'}\lesssim \max\left(\frac{\sqrt{\log(n)}}{n^{1/4}},\frac{h_n^{-\rho'}}{\sqrt{n}}\log(n)\right) \to 0
    \end{align*}
    as $n\to\infty$.

\item We apply Lemma~1 of \cite{GineNickl2008} to show that
\[\cup_{n\ge1}\tilde Q_n=\bigcup_{n\ge1}\left\{\left.x\mapsto\int_{\R}q_t(x-y)K_{h_n}^{(0)}(y)\d y\right|t\in\R\right\}\]
is in the $L^2(\PP)$-closure of $\|K\|_{L^1}$-times the symmetric convex hull of the pregaussian class $Q$. The condition $q_t(\bull-y)\in L^2(\PP)$ is satisfied for all $y\in\R$ since $q_t\in L^2(\R)$ and $ f_Y $ is bounded. $q_t(x-\bull)\in L^1(|\mu_n|)$ is fulfilled owing to $K_{h_n}^{(0)},q_t\in L^2(\R)$. The third condition that $y\mapsto \|q_t(\bull-y)\|_{L^2(\PP)}$ is in $L^1(|\mu_n|)$ holds likewise since $ f_Y $ is bounded and $K_{h_n}^{(0)}\in L^1(\R)$.

\item The $L^2(\PP)$--distance of two functions in $\tilde Q_n$ can be estimated by
\begin{align*}
&     \E\left[(q_t\ast K_h^{(0)}(X)-q_s\ast K_h^{(0)}(X))^2\right]^{1/2}\\
&  =  \left\|\int q_t(\bull-u)K_h^{(0)}(u) - q_s(\bull-u)K_h^{(0)}(u)\d u\right\|_{L^2(\PP)}\\
& \le \int |K_h^{(0)}(u)| \|q_t(\bull-u) - q_s(\bull-u)\|_{L^2(\PP)}\d u\\
& \le \|K_h^{(0)}\|_{L^1}\sup_{|u|\le\xi}  \|q_t(\bull-u) - q_s(\bull-u)\|_{L^2(\PP)}\\
&  =  \|K_h^{(0)}\|_{L^1}\sup_{|u|\le\xi}  \|q_{t+u} - q_{s+u}\|_{L^2(\PP)}.
\end{align*}
As seen in the proof that $Q$ is pregaussian, the covering numbers grow at most polynomially. We take $N$ large enough such that $N\ge2\xi$. Then $s,t>N$ implies $s+u,t+u>N/2$ and $s,t<-N$ implies $s+u,t+u<-N/2$. Since this is only a polynomial change in $N$, the growth of the covering numbers remains at most polynomial. This leads to the entropy bound $H(\tilde Q_n ,L^2(\PP), \eta) \lesssim \log(\eta^{-1})$ for $\eta$ small enough and independent of $n$. We define $\lambda_n(\eta):=\log(\eta^{-1})\eta^2$. The bound in the condition is of the order $\log(n)^{-1/2}n^{1/4}$. As seen before (a) this growth faster than $M_n=C h_n^{-\rho'}$.

\end{enumerate}

\section{Proof of the lower bound}\label{secProofs2}
First we show asymptotic linearity of $\widehat \theta_\zeta$.
\begin{lemma}\label{lemAsympLin}
  Supposing Assumptions~\ref{assFx} and \ref{assEps} and $\zeta\in Z^{\gamma_s,\gamma_c}$ with $\gamma_s>\beta$ and $\gamma_c>(1/2\vee\alpha)+\gamma_s$, the estimator $\widehat\theta_\zeta$ with $h_n=o(n^{-1/(2\alpha+2\gamma_s)})$ is asymptotically linear with influence function $x\mapsto\int\F^{-1}[\phi_\eps^{-1}(-\bull)]\ast\zeta(y)(\delta_x-\PP)(\d y)$ and thus $\widehat\theta_\zeta$ is Gaussian regular.
\end{lemma}
\begin{proof}
  The analysis of the bias of $\widehat\theta$ in Section~\ref{secBias} yields
\begin{align*}
  \widehat \theta=&\theta+\int\F^{-1}[\phi_\eps^{-1}(-\bull)\F K_h]\ast\zeta(y)(\PP_n-\PP)(\d y)+o_P(n^{-1/2})\label{eqRegular}\\
  =&\theta+\int\F^{-1}[\phi_\eps^{-1}(-\bull)]\ast\zeta(y)(\PP_n-\PP)(\d y)\\&+\int\F^{-1}[\phi_\eps^{-1}(-\bull)(\F K_h-1)]\ast\zeta(y)(\PP_n-\PP)(\d y)+o_P(n^{-1/2}).
\end{align*}
Since 
\[\E\left[\left|\int\F^{-1}[\phi_\eps^{-1}(-\bull)]\ast\mathbf\zeta (\d\delta_x-\d\PP)\right|^2\right]\le4\E\left[\int|\F^{-1}[\phi_\eps^{-1}(-\bull)]\ast\mathbf\zeta |^2\d\PP\right]
\]
is finite and $\E[\int\F^{-1}[(\phi_\eps^{-1}(-\bull)]\ast\mathbf\zeta)(\d\delta_x-\d\PP)]=0$ by \eqref{eq2+bound} it suffices to show
\begin{equation}\label{eqOp}
  \int\F^{-1}[\phi_\eps^{-1}(-\bull)(\F K_h-1)]\ast\zeta(y)(\PP_n-\PP)(\d y)=o_P(n^{-1/2}).
\end{equation}
For convenience we write $\psi_h:=\F^{-1}[\phi_\eps^{-1}(-\bull)(\F K_h-1)]\ast\zeta$ and let $\tau>0$. Since $(Y_j)$ are independent and identically distributed, we obtain
\begin{align*}
  &\PP\Big(\big|n^{1/2}\int\psi_h(y)(\PP_n-\PP)(\d y)\big|>\tau\Big)
  \le\tau^{-2}n\E\Big[\big|\int\psi_h(y)(\PP_n-\PP)(\d y)\big|^2\Big]\\
  &\quad=\tau^{-2}n\E\Big[\int\int\psi_h(y)\overline{\psi_h}(z)(\PP_n-\PP)(\d y)(\PP_n-\PP)(\d z)\Big]\\
  &\quad=\tau^{-2}n^{-1}\sum_{j,k=1}^n\E\Big[\int\int\psi_h(y)\overline{\psi_h}(z)(\delta_{Y_j}-\PP)(\d y)(\delta_{Y_k}-\PP)(\d z)\Big]\\
  &\quad=\tau^{-1}\E\Big[\big|\int\psi_h(y)(\delta_{Y_j}-\PP)(\d y)\big|^2\Big]\\
  &\quad\le4\tau^{-1}\int|\psi_h(y)|^2\PP(\d y).
\end{align*}
By uniform integrability of $\psi_h^2$ with respect to $\PP$ by \eqref{eq2+bound} and pointwise convergence $\psi_h\to0$ as $h\to0$ we conclude $\int|\psi_h(y)|^2\PP(\d y)\to0$ and thus \eqref{eqOp}. From asymptotic linearity follows Gaussian regularity by Proposition 2.2.1 of \cite{bicklEtAl1998}.
\end{proof}

Let us now briefly discuss the consequence of Assumption~\ref{assEpsPrim} in terms of Fourier multipliers. Standard calculus yields $|(\phi_\eps^{-1})^{(k)}(u)|\lesssim\langle u\rangle^{\beta-k}$ for $k=0,\dots,(\lfloor\beta\rfloor\vee M)+1$. With the same arguments as in the proof of Lemma~\ref{lemPsiDO}\ref{FourierMultiplier} we deduce that
\begin{align}\label{eqFourierMults}
  (1+iu)^{\beta+k}\phi_\eps^{(k)}(u)\text{ and }(1+iu)^{-\beta+k}(\phi_\eps^{-1})^{(k)}(u)
\end{align}
are Fourier multipliers on $B_{p,q}^s(\R)$ for all $s\in\R,p,q\in[1,\infty]$ and $k=0,\dots\lfloor\beta\rfloor\vee M$.

\subsection{Information bound for smooth $\zeta$}
In this subsection we prove Theorem~\ref{thmCRBound}.\\
\textbf{Step 1:}
To determine the solution of the maximization problem~\eqref{eqCR}, we define $h:=Sg=(g*f_\eps)f_Y^{-1/2}$ with score operator $S$ such that the Fisher information \eqref{FisherInfo} satisfies $\langle \I g,g\rangle=\|h\|^2_{L^2}$. Therefore, we obtain $g=S^{-1}h=\F^{-1}[\varphi_\eps^{-1}]*(\sqrt{f_Y}h)$. Owing to the adjoint equation~\eqref{eqAdjoint}, $\langle g,\zeta\rangle=\int\big(\F^{-1}[\varphi_\eps^{-1}(-\bull)]*\zeta\big)\sqrt{f_Y}h=\langle h,(S^{-1})^\star \zeta\rangle$ holds. Ignoring all restrictions on $g$, the supremum is thus attained at
\begin{equation}\label{eqSInvAdjoint}
  h^*:=(S^{-1})^\star \zeta=(\F^{-1}[\varphi_\eps^{-1}(-\bull)]*\zeta)\sqrt{f_Y}.
\end{equation}
Let us define $\bar\beta:=\lfloor\beta+1/2\rfloor+1$ and $r:=\F^{-1}[
(1+iu)^{-\bar\beta}\varphi_\eps^{-1}(u)]$. Because of Lemma~\ref{lemPsiDO}\ref{lemPsiDOq} we obtain $r\in L^1(\R)\cap L^2(\R)$ and $\F^{-1}[\varphi_\eps^{-1}(u)]=r*(\Id-\D)^{\bar\beta}$. Therefore, the condition $\int g=0$, Fubini's theorem and the fundamental theorem of calculus, provided $(\sqrt{f_Y}h)^{(k)}\in L^1(\R), k=0,\dots\bar\beta$, imply
\begin{align*}
  0=\int\F^{-1}[\varphi_\eps^{-1}]*(\sqrt{f_Y}h)
  &=\int r*\Big(\sqrt{f_Y}h+ \sum_{k=1}^{\bar\beta}\binom{\bar\beta}{k}(-1)^k(\sqrt{f_Y}h)^{(k)}\Big)\\
  &=\int r\Big(\int\sqrt{f_Y}h+\sum_{k=1}^{\bar\beta}\binom{\bar\beta}{k}(-1)^k\int(\sqrt{f_Y}h)^{(k)}\Big).
\end{align*}
For each $k=1,\dots,\bar\beta$ the integrability of $(\sqrt{f_Y}h)^{(l)}, l=k-1,k,$ yields then $\int(\sqrt{f_Y}h)^{(k)}=\lim_{x\to\infty}(\sqrt{f_Y}h)^{(k-1)}(x) -(\sqrt{f_Y}h)^{(k-1)}(-x)=0$ and thus
\begin{align}
  0=\int\F^{-1}[\varphi_\eps^{-1}]*(\sqrt{f_Y}h)
  =\int\sqrt{f_Y}h,\label{eqOrthogonal}
\end{align}
since $\int r=\F r(0)=1$.
Hence, we should project the solution $h^*$ onto the $L^2$-orthogonal space $\operatorname{span}\{\sqrt{f_Y}\}^\perp$:
\begin{align}
  h^{**}&:=h^*-\frac{\langle h^*,\sqrt{f_Y}\rangle}{\|\sqrt{f_Y}\|_{L^2}^2}\sqrt{f_Y}\notag\\
  &=\Big(\F^{-1}[\varphi_\eps^{-1}(-\bull)]*\zeta -\int(\F^{-1}[\varphi_\eps^{-1}(-\bull)]*\zeta)f_Y\Big)\sqrt{f_Y}\nonumber\\
  &=\Big(\F^{-1}[\varphi_\eps^{-1}(-\bull)]*\zeta -\int\zeta f_X\Big)\sqrt{f_Y},\label{eqSg}
\end{align}
where we used $\int(\F^{-1}[\varphi_\eps^{-1}(-\bull)]*\zeta)f_Y=\int\zeta f_X$ by~\eqref{eqAdjoint}. This leads to the candidate for the maximization of \eqref{eqCR} given by
\begin{align*}
  g^*=S^{-1}h^{**}&=S^{-1}(S^{-1})^\star\zeta-\langle\zeta,f_X\rangle S^{-1}\sqrt{f_Y}
  =\I ^{-1}\zeta-\langle\zeta,f_X\rangle f_X\\
  &=\F^{-1}[\varphi_\eps^{-1}]*\Big\{\Big(\F^{-1}[\varphi_\eps^{-1}(-\bull)]*\zeta\Big)f_Y\Big\}-\Big(\int\zeta f_X\Big)f_X
\end{align*}
and \eqref{eqAdjoint} yields $\langle g^{*},\zeta\rangle=\langle \I g^*,g^*\rangle$ and the bound
\begin{equation}\label{eqCRleq}
  \mathcal I_\zeta=\frac{\langle g^*,\zeta\rangle^2}{\langle \I g^*,g^*\rangle}
  =\int\big(\F^{-1}[\phi_\eps^{-1}(-\bull)]*\zeta\big)^2f_Y -\Big(\int\zeta f_X\Big)^2.
\end{equation}
Inequality \eqref{eqBound} holds then by the local version of the H\'ajek--Le~Cam convolution theorem \cite[Thm. 2.3.1]{bicklEtAl1998}. It remains to check the conditions in \eqref{eqCondG},
$(\sqrt{f_Y}h^{**})^{(k)}\in L^1(\R)$ for $k=0,\dots\bar\beta$
and that the three-fold application of the adjoint equality is allowed. The latter will follow from $\sqrt{f_Y}h^{**}, f_Y\in H^{\beta^+}(\R)$ for some $\beta^+>\beta$.\par

\noindent\textbf{Step 2:}
We prove now the integrability of $\sqrt{f_Y}h^{**}=\Big(\F^{-1}[\varphi_\eps^{-1}(-\bull)]*\zeta -\int\zeta f_X\Big)f_Y$ and its derivatives up to order $\bar\beta$ which makes the calculation \eqref{eqOrthogonal} rigorous.\par
For convenience we denote
\[
  \kappa:=\F^{-1}[\phi_\eps^{-1}(-\bull)]*\zeta=r\ast\left(\sum_{k=0}^{\bar\beta}\binom{\bar\beta}{k}(-1)^k\zeta^{(k)}\right).
\]
Owing to Young's inequality together with $r\in L^1(\R)\cap L^2(\R)$ and $\zeta^{(k)}\in L^2(\R)$ for any $k\ge0$, we obtain $\kappa\in C^s(\R)\cap H^s(\R)$ for any $s\ge0$. It suffices to show $f_Y^{(k)}\in L^1(\R)$ for $k=0,\dots,\bar\beta$. Note that by \eqref{eqFourierMults}
\[
  \|(\Id+\D)^kf_\eps\|_{L^1}\lesssim\|\F^{-1}[(1-iu)^k\phi_\eps]\|_{B^0_{1,1}}
  \lesssim\|\F^{-1}[(1-iu)^{k-\beta}]\|_{B^0_{1,1}}
\]
is finite for $\beta>k$ since then $\F^{-1}[(1-iu)^{k-\beta}]=\gamma_{\beta-k,1}\in B_{1,\infty}^{\beta-k}(\R)\subset B_{1,1}^0(\R)$ by the proof of Lemma~\ref{lemPsiDO}\ref{lemPsiDOq}. Recalling that $\beta\notin\Z$, we conclude iteratively $f_\eps^{(k)}\in L^1(\R)$ for $k=0,\dots,\lfloor\beta\rfloor$. Therefore,
\begin{equation*}
  \|f_Y^{(\bar\beta)}\|_{L^1}\le\|f_X^{(\bar\beta-\lfloor\beta\rfloor)}\|_{L^1}\|f_\eps^{(\lfloor\beta\rfloor)}\|_{L^1}<\infty
\end{equation*}
by Assumption~\ref{assFxPrim} and similarly for derivatives of lower order.\par
Moreover, we conclude for $\bar\beta^-\in(\beta+1/2,\bar\beta)$ that 
\[
  f_Y\in B_{1,1}^{\bar\beta^-}(\R)\subset B_{2,1}^{\bar\beta^--1/2}(\R)\subset H^{\beta^+}(\R)
\]
for some $\beta^+>\beta$ by the embeddings \eqref{eqBpEmbed} and \eqref{eqBpEmbed}. Since also $\kappa f_Y\in H^{\beta^+}(\R)$, using $\kappa\in C^{s}(\R)$ for $s>\beta$, we can apply the adjoint equality \eqref{eqAdjoint} in Step 1.\par

\noindent\textbf{Step 3:}
We will show now $\|g^*/f_X\|_\infty<\infty$ which justifies $f_X\pm \tau g^*\ge0$ for some choice of $\tau>0$ small enough.\par
By Step 1 $g^*=\F^{-1}[\varphi_\eps^{-1}]*(\kappa f_Y)-\langle\zeta, f_X\rangle f_X$. For the second term Assumption~\ref{assFxPrim} implies $\|\langle\zeta, f_X\rangle f_X/f_X\|_\infty\le\|\zeta\|_{L^2}\|f_X\|_\infty\|f_X\|_{L^1}<\infty$. Hence, we only need to show $\F^{-1}[\varphi_\eps^{-1}]*(\kappa f_Y)\lesssim f_X$. Using the Besov embedding \eqref{eqCsEmbed}, the Fourier multiplier property of \eqref{eqFourierMults} and the pointwise multiplier property of Besov spaces \eqref{eqPointMult}, we obtain for some $\beta^+\in(\beta,\lfloor\beta\rfloor+1)$
\begin{align*}
  \|\F^{-1}[\varphi_\eps^{-1}]*(\kappa f_Y)\|_\infty
  \lesssim\|\kappa f_Y\|_{B_{\infty,1}^{\beta}}
  \lesssim \|\kappa\|_{B_{\infty,1}^{\beta}}\|f_Y\|_{B_{\infty,1}^{\beta}}
  \lesssim\|\kappa\|_{C^{s}}\|f_Y\|_{C^{\beta^+}}.
\end{align*}
for any $s>\beta$. In Step 2 we have seen that $\kappa\in C^s(\R)$. Moreover,
\begin{align}
  \|f_Y\|_{C^{\beta^+}}
  &=\sum_{k=0}^{\lfloor\beta\rfloor}\|f_Y^{(k)}\|_\infty +\sup_{x\neq y}\Big|\frac{f_Y^{(\lfloor\beta\rfloor)}(x)-f_Y^{(\lfloor\beta\rfloor)}(y)}{(x-y)^{\beta^+-\lfloor\beta\rfloor}}\Big|\notag\\
  &\le\sum_{k=0}^{\lfloor\beta\rfloor}\|f_X\|_\infty \|f_\eps^{(k)}\|_{L^1}+\sup_{x\neq y}\int\frac{|f_X(x-z)-f_X(y-z)|}{|x-y|^{\beta^+-\lfloor\beta\rfloor}}f_\eps^{(\lfloor\beta\rfloor)}(z)\d z\notag\\
  &\le\|f_X\|_\infty\sum_{k=0}^{\lfloor\beta\rfloor} \|f_\eps^{(k)}\|_{L^1}+\|f_X\|_{C^{\beta^+-\lfloor\beta\rfloor}}\|f_\eps^{(\lfloor\beta\rfloor)}\|_{L^1}<\infty,\label{eqFysmooth}
\end{align}
using the Besov embedding $f_X\in W_1^2(\R)\subset B_{1,1}^{\beta^+-\floor\beta+1}\subset C^{\beta^+-\floor\beta}$. Hence, $g^*\in L^\infty(\R)$. Since $f_X$ is a continuous, strictly positive function, we conclude that the quotient $g^*/f_X$ is bounded on every compact subset of $\R$. Therefore, it suffices to estimate the tails. For $|x|$ large enough Assumption~\ref{assFxPrim} implies, using again \eqref{eqFourierMults},
\begin{align*}
  \frac{|\F^{-1}[\varphi_\eps^{-1}]*(\kappa f_Y)(x)|}{f_X(x)}
  &\lesssim \big|x^M\big(\F^{-1}[\varphi_\eps^{-1}]*(\kappa f_Y)\big)(x)\big|\\
  &\le\sum_{k=0}^{M}\tbinom{M}{k}\Big|\F^{-1}\big[(\varphi_\eps^{-1})^{(k)}\F[y^{M-k}\kappa f_Y]\big](x)\Big|\\
  &\lesssim\sum_{k=0}^{M}\|y^{M-k}\kappa f_Y\|_{B_{\infty,1}^{\beta^+-k}}.
\end{align*}
Note that the above calculation shows that $\phi_\eps^{-1}$ is a Fourier multiplier on the weighted Besov space with weight function $\langle x\rangle^M$ \citep[cf. ][Def.~4.2.1/2 and Thm.~5.4.2]{edmundsTriebel1996}. Each term in the above sum can be estimated by
\begin{align*}
  &\|y^{M-k}\kappa\|_{C^s}\|f_Y\|_{C^{\beta^+}}\\
  &\quad=\Big\|\sum_{l=0}^{M-k}\tbinom{M-k}{l}(-1)^l\F^{-1}\big[(\phi_\eps^{-1})^{(l)}(-u)\F[(ix)^{M-k-l}\zeta]\big]\Big\|_{C^s}\|f_Y\|_{C^{\beta^+}},
\end{align*}
where with abuse of notation $\beta^+<\floor\beta+1$ is slightly larger in the last line and $s>\beta^+$. By~\eqref{eqFysmooth} we have $f_Y\in C^{\beta^+}(\R)$. Now, $(ix)^{M-k-l}\zeta\in\mathscr S(\R)$ is again a Schwartz function and thus it suffices to show $\F^{-1}[(\phi_\eps^{-1})^{(k)}(-u)\F\chi]\in C^s(\R)$ for $s>\beta,\chi\in\mathscr S(\R)$ and $k=0,\dots,M$. For $k=0$ this is already done in Step 2. We proceed analogously: for any integer $s\ge0$ we have
\begin{align*}
  &\|\F^{-1}\big[(\phi_\eps^{-1})^{(k)}(-u)\F[\chi]\big]^{(s)}\|_\infty\\
  =&\Big\|\Big(\F^{-1}\big[(1+iu)^{(-\bar\beta+k)\wedge0}(\phi_\eps^{-1})^{(k)}(-u)\big]\ast\big((\Id-\D)^{(\bar\beta-k)\vee0}\chi\big)\Big)^{(s)}\Big\|_\infty\\
  \le&\big\|(1+iu)^{(-\bar\beta+k)\wedge0}(\phi_\eps^{-1})^{(k)}(-u)\big\|_{L^2}\big\|\D^s(\Id-\D)^{(\bar\beta-k)\vee0}\chi\big\|_{L^2}\\
  \lesssim&\big\|\langle u\rangle^{((-\bar\beta+k)\wedge0)+\beta-k}\big\|_{L^2}\big\|\D^s(\Id-\D)^{(\bar\beta-k)\vee0}\chi\big\|_{L^2}.
\end{align*}
Owing to $\bar\beta>\beta+1/2$, the first factor is finite since
\[
  ((-\bar\beta+k)\wedge0)+\beta-k\le
  \begin{cases}
    \beta-\bar\beta<-1/2,\quad&\text{for }\bar\beta\ge k,\\
    \bar\beta-1/2-k<-1/2,\quad&\text{for }\bar\beta< k
  \end{cases}
\]
and the second factor is the $L^2$-norm of a Schwartz function and thus finite, too.

\subsection{Approximation lemma}
To prove convergence of the information bounds it suffices to show that
\begin{align}
  \langle g_n^*,\zeta\rangle&\to\langle (S^\star)^{-1}\zeta,(S^\star)^{-1}\zeta\rangle-\scapro\zeta{f_X}^2\quad\text{and}\label{eqLemKon1}\\
  \langle Sg_n^*,Sg_n^*\rangle=\langle g_n^*,\zeta_n\rangle&\to\langle (S^\star)^{-1}\zeta,(S^\star)^{-1}\zeta\rangle-\scapro\zeta{f_X}^2\label{eqLemKon2}
\end{align}
where we used the equality $\langle g_n^*,\zeta_n\rangle=\langle \I g_n^*,g_n^*\rangle=\scapro{Sg_n^*}{Sg_n^*}$, which holds naturally for the maximizer of the information bound $\mathcal I_{\zeta_n}$. For \eqref{eqLemKon1} we note
\begin{align*}
  \langle g_n^*,\zeta\rangle
  &=\langle \I ^{-1}\zeta_n,\zeta\rangle-\langle\zeta_n, f_X \rangle\langle\zeta,f_X\rangle
  =\langle (S^\star)^{-1}\zeta_n,(S^\star)^{-1}\zeta\rangle-\langle\zeta_n, f_X \rangle\langle\zeta,f_X\rangle
\end{align*}
where the Cauchy--Schwarz inequality yields
\begin{align}
  |\langle f_X, \zeta_n-\zeta\rangle|
  &=|\langle (S^\star)^{-1}(\zeta_n-\zeta),\sqrt{f_Y}\rangle|
  \leq\|(S^\star)^{-1}(\zeta_n-\zeta)\|_{L^2}\to0\label{convergenceTheta}
\end{align}
and
\begin{align*}
  |\langle (S^\star)^{-1}(\zeta_n-\zeta),(S^\star)^{-1}\zeta\rangle|
  &=|\langle (S^\star)^{-1}(\zeta_n-\zeta),(S^\star)^{-1}\zeta\rangle|\\
  &\le\|(S^\star)^{-1}(\zeta_n-\zeta)\|_{L^2}\|(S^\star)^{-1}\zeta\|_{L^2}\to0
\end{align*}
as $n\to\infty$. Analogously follows \eqref{eqLemKon2}, where we use that the assumption of the lemma implies $\langle (S^\star)^{-1}\zeta_n,(S^\star)^{-1}\zeta_n\rangle\to\langle (S^\star)^{-1}\zeta,(S^\star)^{-1}\zeta\rangle$ as $n\to\infty$. The second part of the claim $\theta_{\zeta_n}\to\theta_\zeta$ has already been shown in the estimate~\eqref{convergenceTheta}.

\subsection{Information bound for non-regular $\zeta$}
To prove the efficiency of $\widehat\theta_t$ for $t\in\R$ in Theorem~\ref{thmApprox}, it is suffices by Lemma~\ref{lemApprox} and \eqref{eqSInvAdjoint} to show
\begin{equation}\label{eqApproxZeta}
  \langle (S^\star)^{-1}(\zeta_n-\zeta),(S^\star)^{-1}(\zeta_n-\zeta)\rangle^{1/2}
  =\|\F[\varphi_\eps^{-1}(-\bull)]*(\zeta_n-\zeta)\|_{L^2(\PP)}\to0
\end{equation}
as $n\to\infty$. Using the moment bound \eqref{eq2+bound} replacing $\F K_h$ by 1, we obtain
\[
  \|\F[\varphi_\eps^{-1}(-\bull)]*(\zeta_n-\zeta)\|_{L^2(\PP)}
  \lesssim\|\zeta_n-\zeta\|_{Z^{\beta+\delta,1/2+\beta+\delta}}.
\]
By assumption we have $Z^{\beta+\delta,1/2+\beta+\delta}\subset Z^{\gamma_s,\gamma_c}$ for $\delta$ small enough. Because the space of Schwartz-functions is dense in every Sobolev space $H^s(\R), s\ge0$, $\mathscr S(\R)$ is also dense in $Z^{\gamma_s,\gamma_c}$ and thus the information bound~\eqref{eqBound} holds for all $\zeta\in Z^{\gamma_s,\gamma_c}$. Finally, applying Theorem 25.48 of \cite{vanderVaart1998} and Theorem~\ref{thmPregauss} from above completes the proof of Theorem~\ref{thmApprox}.

\appendix
\section{Appendix: Function spaces}

Let us define the $L^p$-Sobolev space for $p\in(0,\infty)$ and $m\in\N$
\[
W_p^m(\R):=\Big\{f\in L^p(\R)\Big|\sum_{k=0}^m\|f_X^{(k)}\|_{L^p}<\infty\Big\}
\]
In particular, $W^0_p(\R)=L^p(\R)$. Due to the Hilbert space structure, the case $p=2$ is crucial. It can be described equivalently with the notation $\langle u\rangle=(1+u^2)^{1/2}$ by, $\alpha\ge0$,
\[
H^\alpha(\R):=\Big\{f\in L^2(\R)\Big| \|f\|_{H^\alpha}^2:=\int\langle u\rangle^{2\alpha}|\F f(u)|^2\d u<\infty\Big\}
\]
which we call Sobolev space, too. Obviously, $W_2^m(\R)=H^m(\R)$. Also frequently used are the H\"older spaces. Denoting the space of all bounded, continuous functions with values in $\R$ as $C(\R)$ we define, $\alpha\ge0$,
\[
  C^\alpha(\R):=\Big\{f\in C(\R)\Big| \|f\|_{C^\alpha}:=\sum_{k=0}^{\lfloor\alpha\rfloor}\|f^{(l)}\|_\infty+\sup_{x\neq y} \frac{|f^{(\lfloor\alpha\rfloor)}(x)-f^{(\lfloor\alpha\rfloor)}(y)|} {|x-y|^{\alpha-\lfloor\alpha\rfloor}}<\infty\Big\},
\]
where $\lfloor \alpha \rfloor$ denotes the largest integer smaller or equal to $\alpha$.
A unifying approach which contains all function spaces defined so far, is given by Besov spaces \cite[Sect. 2.3.1]{triebel2010} which we will discuss in the sequel. Let $\mathscr S(\R)$ be the Schwartz space of all rapidly decreasing infinitely differentiable functions with values in $\C$ and $\mathscr S'(\R)$ its dual space, that is the space of all tempered distributions. Let $0<\psi\in\mathscr S(\R)$ with $\supp \psi\subset\{x|1/2\le|x|\le2\}$ and $\psi(x)>0$ if $\{x|1/2<|x|<2\}$. Then define $\phi_j(x):=\psi(2^{-j}x)(\sum_{k=-\infty}^\infty\psi(2^{-k}x))^{-1}, j=1,2,\dots,$ and $\phi_0(x):=1-\sum_{j=1}^\infty\phi_j(x)$ such that the sequence $\{\phi_j\}_{j=0}^\infty$ is a smooth resolution of unity. In particular, $\F^{-1}[\phi_j\F f]$ is an entire function for all $f\in\mathscr S'(\R)$. For $s\in\R$ and $p,q\in(0,\infty]$ the Besov spaces are defined by
\[
  B_{p,q}^s:=\Big\{f\in\mathscr S'(\R)\Big|\|f\|_{B^s_{p,q}}:=\Big(\sum_{j=0}^{\infty}2^{sjq}\|\F^{-1}[\phi_j\F f]\|_{L^p}^q\Big)^{1/q}<\infty\Big\}.
\]
We omit the dependence of $\|\bull\|_{B^s_{p,q}}$ to $\psi$ since any function with the above properties defines an equivalent norm. Setting the Besov spaces in relation to the more elementary function spaces, we first note that the Schwartz functions $\mathscr S(\R)$ are dense in every Besov space $B_{p,q}^s$ with $p,q<\infty$ and $H^\alpha(\R)=B_{2,2}^\alpha(\R)$ as well as  $C^\alpha(\R)=B_{\infty,\infty}^\alpha(\R)$, where the latter holds only if $\alpha$ is not an integer \cite[Thms. 2.3.3 and 2.5.7]{triebel2010}. Frequently used are the following continuous embeddings which can be found in \citep[Sect.~2.5.7, Thms. 2.3.2(1), 2.7.1]{triebel2010}: For $p\ge1,m\in\Z$
\begin{align}
  B^m_{p,1}(\R)\subset W^m_p(\R)\subset B^m_{p,\infty}(\R)\quad\text{and}\quad B^0_{\infty,1}(\R)\subset L^\infty(\R)\subset B^0_{\infty,\infty}(\R)\label{eqWpEmbed}
\end{align}
and for $s\ge0$
\begin{align}
  B^s_{\infty,1}(\R)&\subset C^s(\R)\subset B^s_{\infty,\infty}(\R)\label{eqCsEmbed}.
\intertext{Furthermore, for $0<p_0\le p_1\le\infty, q\ge0$ and $-\infty<s_1\le s_0<\infty$}
  B^{s_0}_{p_0,q}(\R)&\subset B^{s_1}_{p_1,q}(\R)\quad\text{if}\quad s_0-\frac{1}{p_0}\ge s_1-\frac{1}{p_1}\label{eqBpEmbed}
\intertext{and for $0<p,q_0,q_1\le\infty$ and $-\infty<s_1< s_0<\infty$}
  B^{s_0}_{p,q_0}(\R)&\subset B^{s_1}_{p,q_1}(\R).\label{eqBqEmbed}
\end{align}
Another important relation is the pointwise multiplier property of Besov spaces \cite[(24) on p. 143]{triebel2010} that is
\begin{equation}\label{eqPointMult}
  \|fg\|_{B^s_{p,q}}\lesssim\|f\|_{B^s_{\infty,q}}\|g\|_{B^s_{p,q}}
\end{equation}
for $s>0$, $1\le p \le \infty$ and $0< q \le\infty$.

The Besov norm of a convolution can be bounded by Lemma~7~(i) in \cite{qui1981}. Let $1\le p,q,r,s\le\infty$, $-\infty<\alpha,\beta<\infty$, $0\le1/u=1/p+1/r-1\le1$, $0\le1/v=1/q+1/s\le1$. For $f\in B^\alpha_{p,q}(\R)$ and $g\in B_{r,s}^\beta(\R)$
\begin{equation}\label{eqConvolutionMult}
\|f\ast g\|_{B^{\alpha+\beta}_{u,v}}\lesssim  \|f\|_{B^\alpha_{p,q}} \|g\|_{ B_{r,s}^\beta}.
\end{equation}

Using for any function $f:\R\to\R$ and $h\in\R$ the difference operators $\Delta_h^1f(x):=f(x+h)-f(x)$ and $(\Delta^l_hf)(x):=\Delta^1_h(\Delta^{l-1}_hf)(x), l\in\N,$ the Besov can be equivalently described by
\[
  \|f\|_{B^s_{pq}}\sim\|f\|_{L^p}+\|f\|_{\dot B_{pq}^s}\quad\text{with}\quad \|f\|_{\dot B_{pq}^s}:=\Big(\int |h|^{-sq-1}\|\Delta_h^Mf\|_{L^p}^q\d h\Big)^{1/q}
\]
for $s>0, p,q\ge1$ and any integer $M>s$ \cite[Thm. 2.5.12]{triebel2010}. The space of all $f\in\mathscr S'(\R)$ for which $\|f\|_{\dot B_{pq}^s}$ is finite is called homogeneous Besov space $\dot B_{pq}^s(\R)$  \cite[Def. 5.1.3/2, Thm. 2.2.3/2]{triebel2010} and thus $B_{pq}^s=L^p(\R)\cap\dot B_{pq}^s(\R)$ for $s>0, p,q\ge1$. Of interest is the relation of homogeneous Besov spaces to functions of bounded $p$-variation. Let $\mathcal{BV}_p(\R)$ denote the space of measurable functions $f:\R\to\R$ such that there is a function $g$ which coincides with $f$ almost everywhere and satisfies
\[
  \sup\Big\{\sum_{i=1}^n|g(x_i)-g(x_{i-1})|^p\Big|-\infty<x_1<\dots<x_n<\infty, n\in\N\Big\}<\infty
\]
and we define $BV_p(\R)$ as the quotient set $\mathcal{BV}_p(\R)$ modulo equality almost everywhere. Then,
\begin{equation}\label{eqPVar}
  \dot B_{p1}^{1/p}(\R)\subset BV_p(\R)\subset \dot B^{1/p}_{p,\infty}(\R), \qquad \text{for } p>1
\end{equation}
by \cite[Thm.~5]{BourdaudEtAl2006}.
For $p=1$ holds by \cite[Lem. 8]{GineNickl2008}
\begin{equation}\label{eqBoundVar}
BV_1(\R)\cap L^1(\R)\subset B^1_{1,\infty}(\R).
\end{equation}

\bibliography{bib} 

\end{document}